\documentclass[11pt,twoside]{article}
\usepackage[utf8]{inputenc}
\usepackage{amsthm,amsmath,amssymb,epsf}

\newcommand{\xx}{{\bf x}}
\newcommand{\yy}{{\bf y}}

\newcommand{\EE}{{\mathbb E}}
\newcommand{\NN}{{\mathbb N}}
\newcommand{\PP}{{\mathbb P}}
\newcommand{\RR}{{\mathbb R}}
\newcommand{\XX}{{\mathbb X}}
\newcommand{\YY}{{\mathbb Y}}

\newcommand{\cC}{{\cal C}}
\newcommand{\cF}{{\cal F}}
\newcommand{\cM}{{\cal M}}
\newcommand{\cT}{{\cal T}}

\newcommand{\1}{{\bf 1}}

\newtheorem{thm}{Theorem}

\newtheorem{cor}{Corollary}
\newtheorem{prop}{Proposition}

\newtheorem{algo}{Algorithm}

\begin{document}
\thispagestyle{empty}

\begin{center}
{\LARGE ABC Shadow algorithm: a tool for statistical analysis of spatial patterns }\\[.5in]
{\large R. S. Stoica$^{1,2}$, A. Philippe$^{3}$, P. Gregori$^4$, J. Mateu$^{4}$}\\[1in]

{\large $^1$Université de Lille, Laboratoire Paul Painlevé, 59655 Villeneuve d'Ascq Cedex, France}\\
{\large $^2$Institut de Mécanique Céleste et Calcul des Ephémérides (IMCCE), Observatoire de Paris, 75014 Paris, France}\\
{\large $^3$ANJA INRIA Rennes Bretagne Atlantique, Université de Nantes, Laboratoire de Mathématiques Jean Leray, 44322 Nantes Cedex 3, France}\\
{\large $^4$Universitat Jaume I de Castell\'{o}n, Instituto Universitario de Matem\'{a}ticas y Aplicaciones de Castell\'{o}n (IMAC), Departamento de Matem\'{a}ticas, Campus Riu Sec, E-12071 Castell\'{o}n, Spain}\\

\end{center}

\begin{verse}
{\footnotesize \noindent {ABSTRACT: }
This paper presents an original ABC algorithm, {\it ABC Shadow}, that can be applied to sample posterior densities that are continuously differentiable. The proposed method uses the ideas given by the auxiliary variable MH of~\cite{MollEtAl06}. The obtained algorithm solves the main condition to be fulfilled by any ABC algorithm, in order to be useful in practice. This condition requires enough samples in the parameter space region, induced by the observed statistics~\cite{Blum10}. The algorithm is tuned on the posterior of a Gaussian model which is entirely known, and then it is applied for the statistical analysis of several spatial patterns. These patterns are issued or assumed to be outcomes of point processes. The considered models are: Strauss, Candy and area-interaction.\\[0.2in]

}
\end{verse}

\section{Introduction}






\noindent
Let us assume that some spatial data set is observed. Typical examples of such sets are digital images, epidemiological or environmental data or catalogues of celestial bodies in astronomy. One typical question related to these examples is the detection and the characterization of the ``hidden'' pattern in the data. Such patterns may be the collection of cells in some biological image, the set of clusters exhibited by a surveyed disease or the filamentary network outlined by the galaxy positions within the observed Universe.\\\

\noindent
Within a probabilistic context, this problem is tackled by assuming that the pattern is the outcome $\yy$ of a stochastic process $\YY$. A possible solution within this context is probabilistic modelling and maximisation. This option requires as a second hypothesis, the knowledge of the model parameters. The model choices may include among others, random fields~\cite{Wink03} or point processes~\cite{Lies00,MollWaag04}.\\

\noindent
Let us consider the observation window to be a finite domain $W \subset \RR^d$ with volume $0 < \nu(W) < \infty$, where $\nu$ is the Lebesgue measure on $\RR^d$. The hidden pattern is supposed to be the realisation $\yy$ of a marked point process $\YY$ on $W \times M$. This means that the hidden pattern is a random set made of a finite number $\kappa$ of random objects $\yy = \{y_1,y_2,\ldots,y_{\kappa}\}$. These objects have a location parameter in $W$ and a characteristic vector or a mark in $M$, that is $y_i=(w_i,m_i) \in W \times M$ for $i=1,\ldots,\kappa$. The marks probability space is $(M, \cM, \nu_{M})$. The pattern probability space is $(\Omega,\cF,\mu)$ with $\Omega$ the configuration space, $\cF$ the associated $\sigma-$algebra and $\mu$ the reference measure.\\

\noindent
The reference measure $\mu$ is given by the unit intensity marked Poisson point process. This process is defined as follows: the number of objects is chosen with respect to a Poisson law of parameter $\nu(W)$, then the objects locations $w_i$ are independently and uniformly chosen in $W$, and finally to each location a mark $m_i$ is attached. The marks are independent of the objects locations and independent of the other marks, and they all follow the same distribution $\nu_M$.\\

\noindent
Due to the independence properties of the process, there is no interaction among the objects in a configuration. Therefore, the outcome of such a process is considered to be a completely random configuration of objects, a completely random pattern. More structured patterns may be obtained by introducing interactions among objects. Such interactions may be specified by means of a probability density $p$ with respect to the reference measure $\mu$. The Gibbsian modelling framework allows to write such a probability density as 
\begin{equation}
p(\yy|\theta) = \frac{\exp[-U(\yy|\theta)]}{c(\theta)}
\label{modelPattern}
\end{equation}
with $U : \Omega \rightarrow \RR^{+}$ the energy function, $\theta$ the model parameters and $c(\theta)$ the normalising constant.\\

\noindent
Being in the possession of a model~\eqref{modelPattern}, the hidden pattern estimator is given by
\begin{equation}
\widehat{\yy} = \arg\max_{\yy \in \Omega}\{p(\yy|\theta)\} 
= \arg\min_{\yy \in \Omega}\{U(\yy|\theta)\}.
\label{patternEstimator}
\end{equation}
The optimisation in~\eqref{patternEstimator} is done using a simulated annealing algorithm based on a Metropolis-Hastings (MH) or spatial birth-and-death dynamics~\cite{Lies94,StoiGregMate05}. It is important to notice that these algorithms do not require the computation of the normalising constant $c(\theta)$ which is not always available in analytical closed form.\\

\noindent
The dual formulation of the pattern detection question is the parameter estimation problem. Let us now  consider that an object pattern $\yy$ is observed in $W$. The observed pattern is supposed to be the realisation of a marked point process given by the probability density $p(\yy|\theta)$. Let $p(\theta|\yy)$ be the conditional distribution of the model parameters or the posterior law
\begin{equation} 
p(\theta|\yy) = \frac{\exp[-U(\yy|\theta)]p(\theta)}{Z(\yy)c(\theta)},
\label{posteriorGibbs}
\end{equation}
$p(\theta)$ the prior density for the model parameters and $Z(\yy)$ the normalising constant. The posterior law is defined on the parameter space $\Theta$. For simplicity, the parameter space is considered to be a compact region in $\RR^{r}$ with $r$ the size of the parameter vector. The parameter space is endowed with its Borel algebra $\cT_{\Theta}$.\\

\noindent
Straightforward sampling of~\eqref{posteriorGibbs}, for instance using a classical MH dynamics, cannot be always done, since it may require the evaluation of the ratio $c(\theta)/c(\psi)$ for $(\theta,\psi) \in \Theta^2$. The theoretical solution to this problem is given by~\cite{MollEtAl06}. The authors propose an auxiliary variable MH sampler that, by a clever and elegant choice of the proposal and the auxiliary variable densities, avoids the computation of the normalising constants. Nevertheless, as the authors themselves explain in detail, it is difficult to show how these choices prevent the simulated chain from poor mixing. Approximate solutions for posterior sampling are provided by Approximate Bayesian Computation (ABC) methods~\cite{BeauEtAl09,GrelEtAl09,MariEtAl12,AtchLartRobe13}. These methods are very attractive from a practical point of view. The theoretical developments of~\cite{Blum10,BiauEtAl15} established nice connections of the ABC based inference with kernel and $k-$nearest neighbour estimation. The bias and the variance for the posterior distribution estimate is given by~\cite{Blum10}. Despite its simplicity, the main criticism against ABC methods is the large number of generated samples that cannot be used. In fact, in order to get reliable numerical results, an ABC method should produce enough samples in the parametric region induced by the posterior sampling.\\

\noindent
The aim of this paper is to solve this drawback for posterior densities that are continuously differentiable. This a strong but realistic assumption. With this goal in mind, we embed the auxiliary sampling ideas of~\cite{MollEtAl06} in the ABC framework given by~\cite{Blum10,BiauEtAl15}. The paper continues with the presentation of the auxiliary variable method of~\cite{MollEtAl06}. Next the ABC methods from~\cite{Blum10,BiauEtAl15} are described. The principles of our posterior sampling method are shown in Section~3. The application results are analysed in Section~4. The method is first tuned on a Gaussian model for which the posterior density is entirely known. Then the proposed algorithm is applied to the statistical analysis of spatial patterns issued or assumed to be point processes outcomes. The considered models are: the Strauss, the Candy and the area-interaction processes. At our best knowledge, it is the first time an ABC based statistical analysis is done for such models. The paper ends with conclusions and perspectives.\\

\section{Sampling posterior distributions}
\subsection{Theoretical solution: auxiliary variable MH algorithm}
\noindent
The goal is to sample the posterior
\begin{equation*}
p(\theta|\yy) \propto p(\yy|\theta)p(\theta).
\end{equation*}

\noindent
The authors idea in~\cite{MollEtAl06} is to use an auxiliary variable $\XX$ with probability density $a(\xx|\theta,\yy)$. Consequently, the joint distribution becomes
\begin{equation}
\pi(\theta,\xx|\yy) = a(\xx|\theta,\yy)p(\theta|\yy) = a(\xx|\theta,\yy)\frac{\exp[-U(\yy|\theta)]p(\theta)}{Z(\yy)c(\theta)}.
\label{posteriorAuxLaw}
\end{equation}

\noindent
Sampling from~\eqref{posteriorAuxLaw} can be done using a MH algorithm that uses a proposal made of two components, such as:
\begin{equation*}
q((\theta,\xx) \rightarrow (\theta^{\prime},\xx^{\prime})) = q_{1}(\theta^{\prime}|\theta,\xx)q_{2}(\xx^{\prime}|\theta^{\prime},\theta,\xx).
\end{equation*}
The first component $q_{1}$ may be a very simple probability law and it may not depend on $\xx$. Let us consider, for instance, the symmetric distribution 
\begin{equation}
q_{1}(\theta^{\prime}|\theta) = q_{1}(|\theta^{\prime} - \theta|).
\label{danishProposal1}
\end{equation}
The second component $q_{2}$ is chosen independently on $(\theta,\xx)$ as follows 
\begin{equation}
q_{2}(\xx^{\prime}|\theta^{\prime},\theta,\xx) = q_{2}(\xx^{\prime}|\theta^{\prime}) = \frac{\exp[-U(\xx^{\prime}|\theta^{\prime})]}{c(\theta^{\prime})}.
\label{danishProposal2}
\end{equation}

\noindent
The MH ratio $H_D((\theta,\xx) \rightarrow (\theta^{\prime},\xx^{\prime}))$ of the induced dynamics is, by definition,
\begin{equation}
H_D((\theta,\xx) \rightarrow (\theta^{\prime},\xx^{\prime})) = \frac{\pi(\theta^{\prime},\xx^{\prime}|\yy)}{\pi(\theta,\xx|\yy)} \times \frac{q((\theta^{\prime},\xx^{\prime}) \rightarrow (\theta,\xx))}{q((\theta,\xx) \rightarrow (\theta^{\prime},\xx^{\prime}))}.
\label{danishMH}
\end{equation}

\noindent
Plugging the joint density~\eqref{posteriorAuxLaw} and the proposals~\eqref{danishProposal1} and~\eqref{danishProposal2} into~\eqref{danishMH}, the MH ratio becomes
\begin{equation*}
\frac
{a(\xx^{\prime}|\theta^{\prime},\yy)\exp[-U(\yy|\theta^{\prime})]p(\theta^{\prime})}
{a(\xx|\theta,\yy)\exp[-U(\yy|\theta)]p(\theta)}
\times
\frac
{\exp[-U(\xx|\theta)]}
{\exp[-U(\xx^{\prime}|\theta^{\prime})]},
\end{equation*}
which does not depend on the normalising constant $c$.\\


\noindent
For exponential family models, the MH ratio has a simpler form. In this case, the energy functions are of the form $U(\yy|\theta) = \langle t(\yy),\theta \rangle$, where $\langle,\rangle$ is the scalar product of the sufficient statistics vector of the model $t(\yy)$ and the parameter vector $\theta$, so the MH ratio is
\begin{equation*}
\frac{a(\xx^{\prime}|\theta^{\prime},\yy) p(\theta^{\prime}) \exp[ - \langle t(\yy), \theta^{\prime} \rangle] \exp[ - \langle t(\xx), \theta \rangle]}
{a(\xx|\theta,\yy) p(\theta) \exp[ - \langle t(\yy), \theta \rangle] \exp[ - \langle t(\xx^{\prime}), \theta^{\prime} \rangle]}
\end{equation*}

\noindent
Under these considerations, the auxiliary variable MH sampler works as follows:

\begin{algo}
Assume the observed pattern is $\yy$ and the current state is $(\theta_k,\xx_{k})$.
\begin{enumerate}
\item{Generate a new candidate $\theta^{\prime}$ using the proposal $q_{1}(\theta^{\prime}|\theta)$.}
\item{Generate a new candidate $\xx^{\prime}$ using the proposal $q_{2}(\xx^{\prime}|\theta^{\prime}) = \frac{\exp[-U(\xx^{\prime}|\theta^{\prime})]}{c(\theta^{\prime})}$.}
\item{Compute the MH ratio $H_D((\theta,\xx) \rightarrow (\theta^{\prime},\xx^{\prime}))$ using~\eqref{danishMH} or one of its simplified forms previously detailed.}
\item{The new state $(\theta_{k+1},\xx_{k+1})=(\theta^{\prime},\xx^{\prime})$ is accepted with probability
\begin{equation*}
\alpha = \min\left\{1,H_D((\theta,\xx) \rightarrow (\theta^{\prime},\xx^{\prime}))\right\},
\end{equation*}
otherwise, the chain remains in the same initial state, that is 
\begin{equation*}
(\theta_{k+1},\xx_{k+1})=(\theta_k,\xx_k).
\end{equation*}
}
\item{If more samples are needed, re-iterate the whole procedure.}
\end{enumerate}
\label{auxiliaryVariableAlgorithm}
\end{algo}

\noindent
The outputs of the Algorithm~\ref{auxiliaryVariableAlgorithm} converge in law towards the unique equilibrium distribution $\pi(\theta,\xx|\yy)$~\cite{MollEtAl06}. 
This result holds if sampling from the proposal $q_{2}(\xx^{\prime}|\theta^{\prime})$ is done exactly. Despite the fact that the convergence time of the exact simulation methods for marked point processes may highly depend on model parameters~\cite{LiesStoi06}, we do not consider this characteristic as a drawback of the auxiliary variable method, since this method is mathematically correct. We hope that future advances in computer science and algorithmics will be able to propose a practical implementation of this solution. The difficulty in a straightforward use of the auxiliary variable method is the freedom of choice of the auxiliary variable density. As the authors themselves showed in~\cite{MollEtAl06}, there is no clear choice for it, such that good mixing properties of the simulated chain are guaranteed. A possible strategy, that may be fair enough, is to set $a(\xx|\theta,\yy) = \frac{\exp[-U(\xx|\widehat{\theta}(\yy))]}{c(\widehat{\theta}(\yy))}$ where $\widehat{\theta}(\yy)$ is the pseudo-likelihood estimate of the model parameters based on the observation of $\yy$.\\

\subsection{Approximate solution: ABC methods}
\noindent
ABC (Approximate Bayesian Computation) is the generic name for numerical simulation methods allowing statistical inference based on an approximate sampling from the posterior distribution $p(\theta|\yy)$~\cite{AtchLartRobe13,BeauEtAl09,GrelEtAl09,MariEtAl12}. Its general idea is described by the following algorithm.\\

\begin{algo}
Assume the observed pattern is $\yy$, fix a tolerance threshold $\epsilon$ and an integer value $n$.
\begin{enumerate}
\item For $i=1$ to $n$ do
\begin{itemize}
\item Generate $\theta_i$ according to $p(\theta)$.
\item Generate $\xx_{i}$ according to the probability density 
\begin{equation*}
p(\xx|\theta_i) = \frac{\exp[-U(\xx|\theta_i)]}{c(\theta_i)}
\end{equation*}
\end{itemize}
\item Return all the $\theta_i$'s such that the distance between the statistics of the observation and those of the simulated pattern is small, that is 
\begin{equation*}
d(t(\yy),t(\xx_i)) \leq \epsilon
\end{equation*}
\end{enumerate}
\label{abcAlgorithm}
\end{algo}

\noindent
In the general case, the outputs of the Algorithm~\ref{abcAlgorithm} are distributed according to the law $\pi(\theta|d(t(\yy),t(\xx)) \leq \epsilon)$. The author in~\cite{Blum10} established a link between kernel estimation and inference based on this algorithm. In this paper, estimates of the bias and the variance of the posterior density estimate are given. The choices influencing the algorithm outputs are related to the statistics vector, the distance definition and the tolerance threshold. For exponential family models the sufficient statistics vector is the natural choice.\\

\noindent
The authors in~\cite{BiauEtAl15} considered a slightly different version of the previous algorithm.

\begin{algo}
Assume the observed pattern is $\yy$ and two integer values $k_n$ and $n$ such that $1 \leq k_n \leq n$.
\begin{enumerate}
\item For $i=1$ to $n$ do
\begin{itemize}
\item Generate $\theta_i$ according to $p(\theta)$.
\item Generate $\xx_{i}$ according to the probability density 
\begin{equation*}
p(\xx|\theta_i) = \frac{\exp[-U(\xx|\theta_i)]}{c(\theta_i)}
\end{equation*}
\end{itemize}
\item Return all the $\theta_i$'s such that $t(\xx_i)$ is among the $k_n$-nearest neighbours of $t(\yy)$.
\end{enumerate}
\label{abcAlgorithm-NN}
\end{algo}

\noindent
The Algorithms~\ref{abcAlgorithm-NN} and~\ref{abcAlgorithm} are similar. The same requirements hold for both algorithms. But for the nearest neighbour algorithm no tolerance threshold is needed. Instead, the $k_n$ parameter should be fixed. Intuitively, the smaller $\epsilon$ is or the higher $k_n$ is (such that $k_n/n \rightarrow 0$ with $n \rightarrow \infty$), the closer the algorithms outputs are to the posterior density. In both cases, the posterior density is estimated using a kernel of a given bandwidth, over the selected samples. Clearly, the bandwidth of the kernel influences the quality of the estimation. This is not a drawback in our opinion, whenever enough points are sampled close to the observed statistics. So, as indicated in~\cite{Blum10}, the key point in building an ABC procedure is to ensure that the sampling mechanisms of the algorithm do not produce too many samples that are sparse around the observed statistics.

\section{ABC Shadow algorithm}
The ABC method we propose follows directly the recommendation of~\cite{Blum10}, to sample close to the posterior. Our solution is inspired by the auxiliary variable method~\cite{MollEtAl06}. This goal is achieved by building two Markov chains, the ideal and the shadow chains. The ideal chain is a theoretical chain that cannot be used in practice, but its equilibrium distribution is the posterior law we want to sample from. The shadow is a chain that can be practically simulated, and it follows as close as desired the ideal chain during a fixed number of steps.

\subsection{Ideal chain}
In theory, Markov chain Monte Carlo algorithms may be used for sampling $p(\theta | \yy)$. For instance, let us consider the general MH algorithm. Assuming the system is in the state $\theta$, this algorithm first chooses a new value $\psi$ according to a proposal density $q(\theta \rightarrow \psi)$. The value $\psi$ is then accepted with probability $\alpha_{i}(\theta \rightarrow \psi)$ given by
\begin{equation}
\alpha_{i}(\theta \rightarrow
\psi)=\min\left\{1,\frac{p(\psi|\yy)}{p(\theta|\yy)}\frac{q(\psi
\rightarrow \theta)}{q(\theta \rightarrow \psi)}\right\}.
\label{acceptance_probability_ideal}
\end{equation}

\noindent 
The transition kernel of the Markov chain simulated by this algorithm is
\begin{eqnarray*}
P_i(\theta,A) & = &\int_{A}\alpha_{i}(\theta \rightarrow
\psi)q(\theta \rightarrow \psi)\1\{\psi \in A\}d\psi\\
& +  &\1\{\theta \in A\}
\left[1 - \int_{A}\alpha_{i}(\theta \rightarrow
\psi)q(\theta \rightarrow \psi) d\psi \right]
\label{transition_kernel}
\end{eqnarray*}
with $A \in \cT_\Theta$.\\

\noindent
The conditions that the proposal density $q(\theta \rightarrow \psi)$ has to meet, so that the simulated Markov chain has a unique equilibrium distribution
\begin{equation*}
\pi(A)=\int_{A}p(\theta|\yy)d\nu(\theta),
\end{equation*}
\noindent
are rather mild~\cite{Tier94}. Furthermore, the simulated chain is uniformly ergodic, that is, there is a positive constant $M$ and a positive constant $\rho < 1$ such that
\begin{equation*}
\sup_{\theta \in \Theta}\parallel P_{i}^{n}(\theta,\cdot) - \pi(\cdot) \parallel \leq M\rho^{n}, \quad n \in \NN.
\end{equation*}

\noindent
For a fixed $\Delta > 0$, a parameter value $\nu \in \Theta$ and a realisation $\xx$ of the model $p(\cdot|\nu)$ given by~\eqref{modelPattern}, let us consider the proposal density
\begin{eqnarray}
q(\theta \rightarrow \psi) = q_\Delta (\theta \rightarrow \psi |
\xx) = \frac{f(\xx |\psi)/c(\psi)}{I(\theta, \Delta,
\xx)}\1_{b(\theta, \Delta/2)}\{\psi\} \label{proposal_density}
\label{idealProposal}
\end{eqnarray}
with $f(\xx|\psi) = \exp[-U(\xx|\psi)]$. Here $\1_{b(\theta, \Delta/2)}\{\cdot\}$ is the indicator function over $b(\theta, \Delta/2)$, which is the ball of centre $\theta$ and radius $\Delta/2$. Finally, $I(\theta, \Delta, \xx)$ is the quantity given by the integral 
\begin{equation*}
I(\theta, \Delta, \xx) = \int_{b(\theta, \Delta/2)} f(\xx|\phi)/c(\phi) \, d \phi.
\end{equation*}

\noindent
This choice for $q(\theta \rightarrow \psi)$ guarantees the convergence of the chain towards $\pi$ and avoids the evaluation of the normalising constant ratio $c(\theta)/c(\psi)$ in~\eqref{acceptance_probability_ideal}. We call the chain induced by these proposals the {\it ideal chain}. Nevertheless, the proposal~\eqref{idealProposal} requires the computation of integrals such as $I(\theta, \Delta, \xx)$, and this is as difficult as the computation of the normalising constant ratio. Still, this construction allows a natural approximation of the ideal chain: the shadow chain.\\

\subsection{Shadow chain}
\noindent
In order to ease the reading of this section, the proofs of the theorems and propositions are given in the Appendix at the end of the paper.\\

\noindent
Let $V_\Delta$ denote the volume of the ball $b(\theta, \Delta/2)$, and let 
\begin{equation}
U_\Delta (\theta \to \psi) = \frac{1}{V_\Delta} \1_{b(\theta, \Delta/2)}\{\psi\}
\label{uniformProposal}
\end{equation}
be the uniform probability density over the ball $b(\theta, \Delta/2)$. The probability density $p(\xx|\phi)$ defined by~\eqref{modelPattern} is assumed to be a continuously differentiable function in $\phi$, $p(\xx|\cdot) \in \cC^{1}(\Theta)$.\\

\begin{thm} 
\label{first_theorem}
Let $\xx$ be a point in $\Omega$ such that the
function $p(\xx | \phi)$ is strictly positive and continuous with
respect to $\phi$, then we have that:
\renewcommand{\labelenumi}{$($\roman{enumi}$)$}
\begin{enumerate}
\item The probability distributions given by the proposal densities $q_\Delta( \theta \to \cdot )$ and $U_\Delta ( \theta \to \cdot )$  are uniformly as close as desired in $\theta$ as $\Delta$ approaches $0$. That is, for any fixed $\theta \in \Theta$ and $A \in \cT_{\Theta}$, we have
\begin{equation*}
\lim_{\Delta \to 0_+} \int_A | q_{\Delta}( \theta \to \psi ) -
U_{\Delta}( \theta \to \psi ) | d \psi = 0.
\end{equation*}
\item
For any fixed $\theta \in \Theta$, the quotient functions
$\frac{{q}_\Delta( \theta \to \cdot )}{{q}_\Delta( \cdot \to \theta
)}$ and $\frac{
\frac{f(\xx|\cdot)}{c(\cdot)}\1_{b(\theta,\Delta/2)}(\cdot) }{
\frac{f(\xx|\theta)}{c(\theta)}\1_{b(\cdot,\Delta/2)}(\theta) }$ are
uniformly as close as desired in $\theta$ as $\Delta$ approaches
$0$. That is, for any fixed $\theta \in \Theta$, we have
\begin{equation*}
\lim_{\Delta \rightarrow 0_{+}} \sup_{\psi \in \Theta} \left|
\frac{q_\Delta (\theta \rightarrow \psi | \xx)}{q_\Delta (\psi
\rightarrow \theta | \xx)} - \frac{
\frac{f(\xx|\psi)}{c(\psi)}\1_{b(\theta,\Delta/2)}(\psi) }{
\frac{f(\xx|\theta)}{c(\theta)}\1_{b(\psi,\Delta/2)}(\theta) }
\right| = 0
\end{equation*}
uniformly in $\theta \in \Theta$.
\end{enumerate} 
Moreover, if $p(\xx|\cdot) \in \cC^{1}(\Theta)$, then rates of convergence in (i) and (ii) can be provided.
\end{thm}

\noindent
This result leads to the construction of a new Markov chain that approximates the behaviour of the ideal chain for small values of $\Delta$. The first part of the theorem allows the use of the uniform law $U_{\Delta}(\theta \rightarrow \psi)$ instead of $q_{\Delta}(\theta \rightarrow \psi)$ for proposing new values. The second part of the result enables the approximation of the ideal acceptance ratio. The resulting new chain is called the {\it shadow chain}. Assuming the chain in a state $\theta$, a new value $\psi$ is chosen uniformly in the ball $b(\theta,\Delta/2)$. The state $\psi$ is accepted with probability
\begin{equation}
\alpha_{s}(\theta \rightarrow \psi) =
\min\left\{1,\frac{p(\psi|\yy)}{p(\theta|\yy)}\times\frac{f(\xx | \theta)c(\psi)\1_{b(\psi,\Delta/2)}\{\theta\}}
{f (\xx | \psi)c(\theta)\1_{b(\theta,\Delta/2)}\{\psi\}}\right\},
\label{acceptance_probability_shadow}
\end{equation}
otherwise the chain remains in the initial state.\\

\noindent
By construction, the shadow chain is irreducible and aperiodic~(\cite{MeynTwee09}). Yet, we do not have any knowledge about the existence and the uniqueness of its equilibrium distribution. The following corollary is a direct consequence of Theorem~\ref{first_theorem}.\\

\begin{cor}
\label{first_corollary}
The acceptance probabilities of the ideal and shadow chains given by~\eqref{acceptance_probability_ideal} and~\eqref{acceptance_probability_shadow} respectively, are uniformly as closed as desired whenever $\Delta$ approaches $0$.
\end{cor}

\noindent
Next, we show that, for a given $\xx$, it is possible to choose $\Delta$ such that during $n$ steps, the ideal and shadow chains evolve as close as desired.

\begin{prop}
Let $P_{i}$ and $P_{s}$ be the transition kernels for
the ideal and the shadow Markov chains using a general $\Delta > 0$ and a configuration $\xx \in \Omega$ as in Theorem~\ref{first_theorem},
respectively. Then for every $\epsilon >0$ and every $n \in \NN$,
there exists $\Delta_0 = \Delta_0(\epsilon, n) >0$ such that for
every $\Delta \leq \Delta_0$ we have $| P_{i}^{(n)}(\theta,
A) - P_{s}^{(n)}(\theta, A) | < \epsilon$ uniformly in
$\theta \in \Theta$ and $A \in \cT_{\Theta}$. If $p(\xx|\theta) \in \cC^1(\Theta)$, then a description of $\Delta_0(\epsilon, n)$ can be provided.
\label{first_proposition}
\end{prop}

\noindent
The previous results allow the construction of the following ABC algorithm. Its outputs are approximate samples from the posterior $p(\theta|\yy)$.\\

\begin{algo}
\label{abcShadow}
{\bf ABC Shadow~:} fix $\Delta$ and $n$. Assume the observed pattern is $\yy$ and the current state is $\theta_0$.\\

\begin{enumerate}
\item Generate $\xx$ according to $p(\xx|\theta_0)$.
\item For $k=1$ to $n$ do\\
\begin{itemize}  
\item Generate a new candidate $\psi$ following $U_\Delta(\theta_{k-1} \to \psi)$ defined by~\eqref{uniformProposal}.
\item The new state $\theta_{k} = \psi$ is accepted with probability $\alpha_{s}(\theta_{k-1} \rightarrow \psi)$ given by~\eqref{acceptance_probability_shadow}, otherwise $\theta_{k} = \theta_{k-1}$.
\end{itemize}
\item Return $\theta_n$
\item If another sample is needed, go to step $1$  with $\theta_0 = \theta_n$.
\end{enumerate}
\end{algo}

\noindent
If $n \rightarrow \infty$ the Algorithm~\ref{abcShadow} may not follow closely an ideal chain started in $\theta_0$. In this case, the shadow chain drifts away from the stationary regime of the ideal chain. Nevertheless, for a fixed $n$ and $\Delta \rightarrow 0$, the ideal chain approaches the equilibrium regime, equally fast from any initial point. If $\epsilon$, the distance after $n$ steps between the ideal and the shadow chain is small enough, the shadow chain will also get closer to the equilibrium regime. The triangle inequality gives a bound for the distance after $n$ steps, between the shadow transition kernel and the equilibrium regime
\begin{equation*}
\| P_s^{(n)}(\theta,\cdot) - \pi(\cdot) \| \leq M(\xx,\Delta)\rho^n + \epsilon.
\end{equation*}

\noindent
One very important element here is the auxiliary variable $\xx$. Refreshing it allows to re-start the algorithm for another $n$ steps more, and by this, to obtain new samples of the approximate distribution. Comparing with classical ABC, no useless samples are produced, in the sense that the outputs of the Shadow algorithm tend to have a distribution ``closer'' to the true posterior while $\xx$ is renewed. In the next section, a simulation study is done, in order to analyse the Shadow algorithm performances.\\

\section{Applications}
\subsection{Method analysis: posterior approximation for a Gausian model}
\noindent
The Algorithm~\ref{abcShadow} can be applied to exponential family models. Here, the method is tested on the posterior of a Gaussian model. In this case, the posterior distribution has an analytical expression allowing its simulation using a MH dynamics. The comparison of the ABC Shadow algorithm with the classical MH dynamics allows the tuning of the algorithm and its numerical analysis in terms of mixing, parameter sensibility and initial conditions dependence.\\

\noindent
The posterior distribution of a Gaussian model with mean $\theta_1$ and variance $\theta_2$ is
\begin{equation}
p(\theta_1,\theta_2|\yy ) \propto \frac{\exp\left(\frac{\theta_1}{\theta_2}t_{1}(\yy) - \frac{t_{2}(\yy)}{2\theta_2}\right)}{c(\theta_1,\theta_2)}p(\theta_1,\theta_2),
\label{posteriorNormal}
\end{equation}
with $t(\yy)=(t_{1}(\yy),t_{2}(\yy))$ the observed sufficient statistics vector. If a sample of size $m$ is observed, that is the sequence $\yy = \{\YY_{1}(\omega),\YY_{2}(\omega),\ldots,\YY_{m}(\omega)\}$ is given by $m$ independent random variables following a Gaussian distribution with parameters $\theta_1$ and $\theta_2$, then the sufficient statistics vector is 
\begin{equation*}
t(\yy) = \left ( \sum_{i=1}^{m}\YY_{i}(\omega), \sum_{i=1}^{m}\YY_{i}^{2}(\omega)\right).
\end{equation*}
 
\noindent
For our experiment, we simulated $m=1000$ independent identical Normal random variables with parameters $\theta=(\mu,\sigma^2)=(2,9)$. The observed sufficient statistics vector is $t(\yy) = (1765.45,12145.83)$ and $p(\theta)$ the uniform distribution over the interval $[-100,100] \times [0,200]$. The MH algorithm parameters implemented to sample from~\eqref{posteriorNormal} used uniform proposals of width $(0.5,0.5)$ around the current values. Finally, this dynamics was run for $125 \times 10^5$ iterations and samples were kept every $12500$ steps. This gave an amount of $1000$ samples. For the ABC Shadow algorithm, the $\Delta$ parameter was set to $(0.005,0.025)$ and $n=500$. The algorithm was run $25 \times 10^3$ times. Samples were kept every $25$ repetitions of the ABC procedure. This gave an amount of $1000$ samples as for the MH dynamics.\\

\begin{figure}[!htbp]
\begin{center}
\begin{tabular}{cc}
\epsfxsize=6cm \epsffile{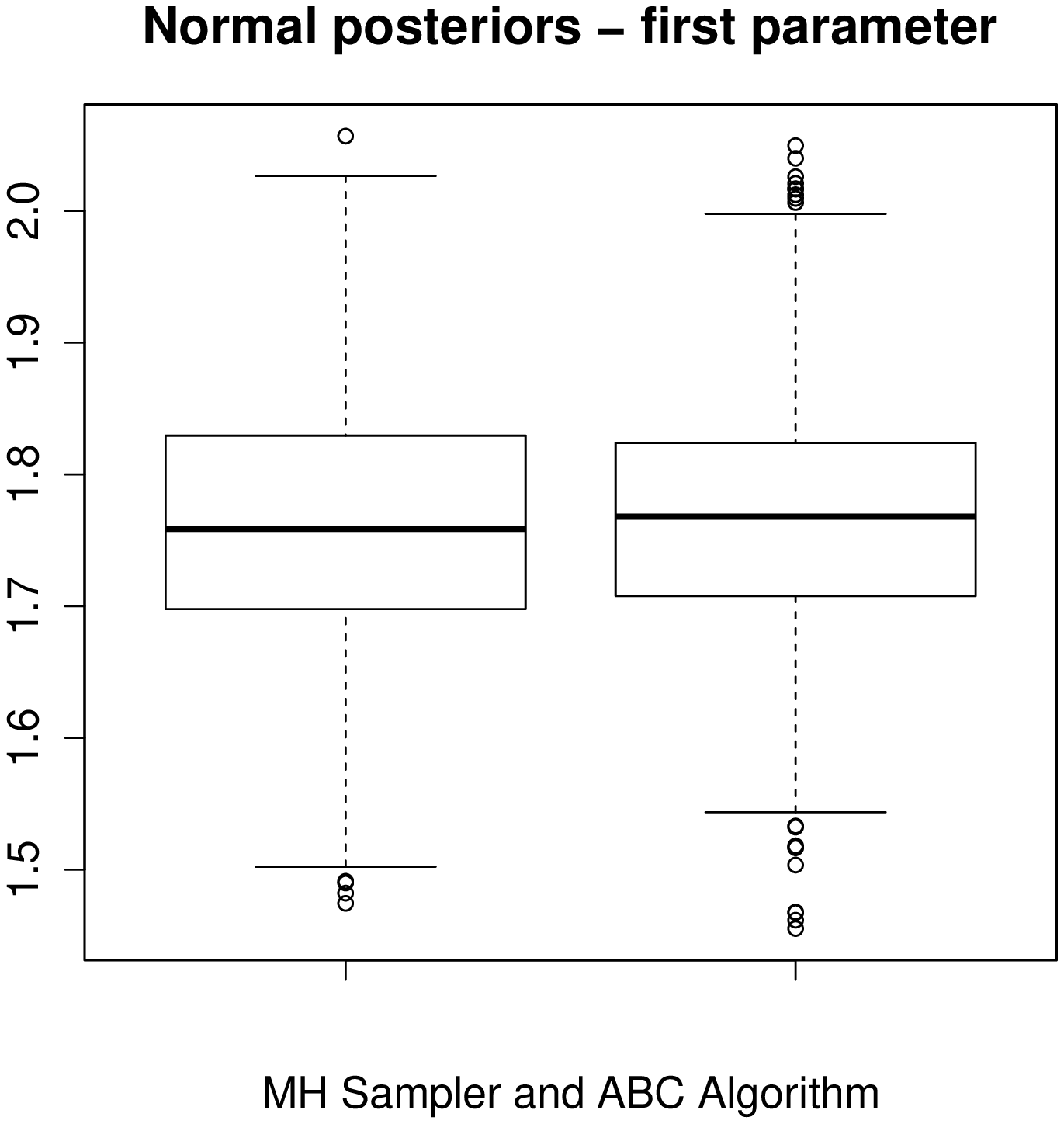} & \epsfxsize=6cm \epsffile{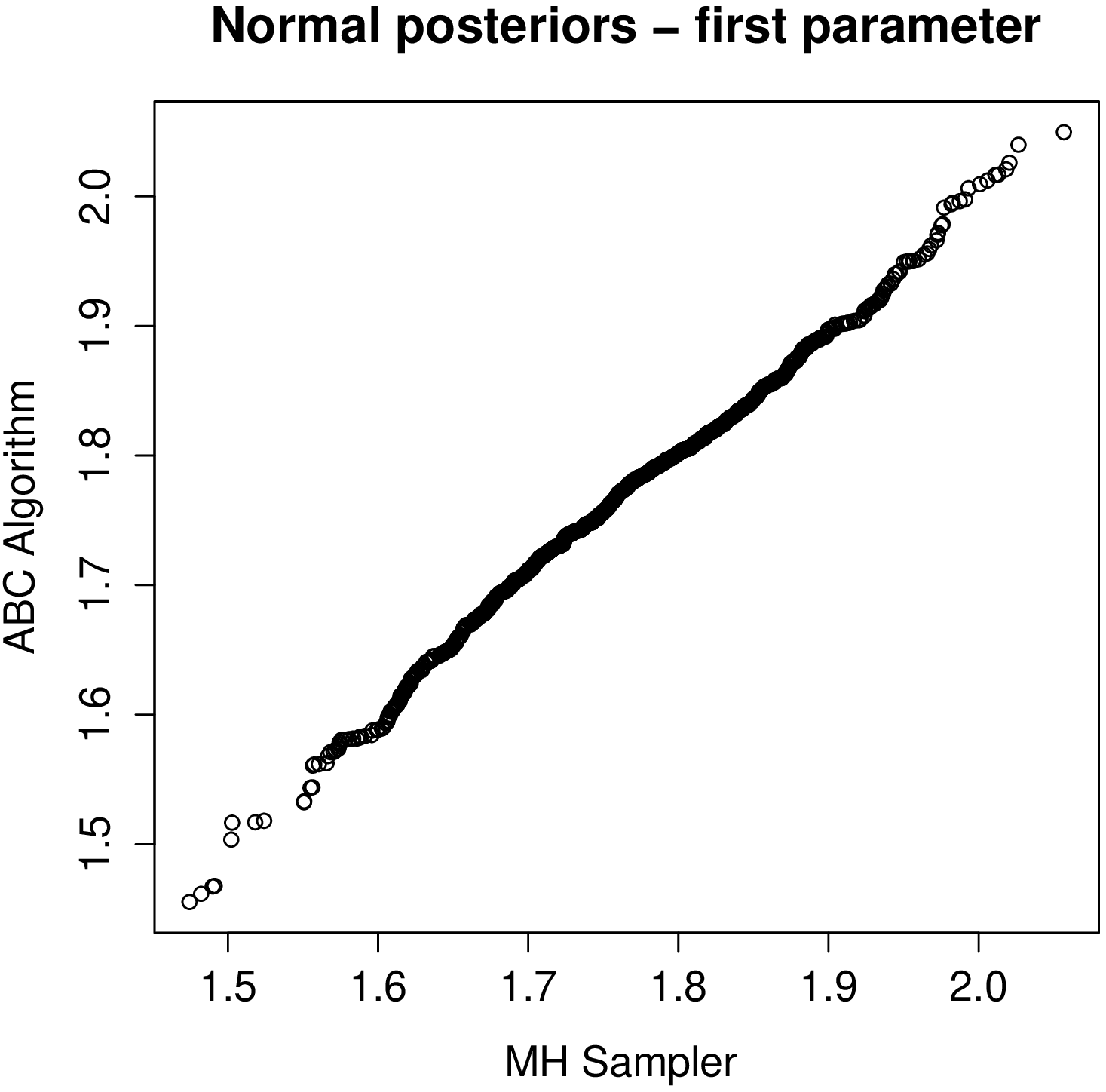} \\
\epsfxsize=6cm \epsffile{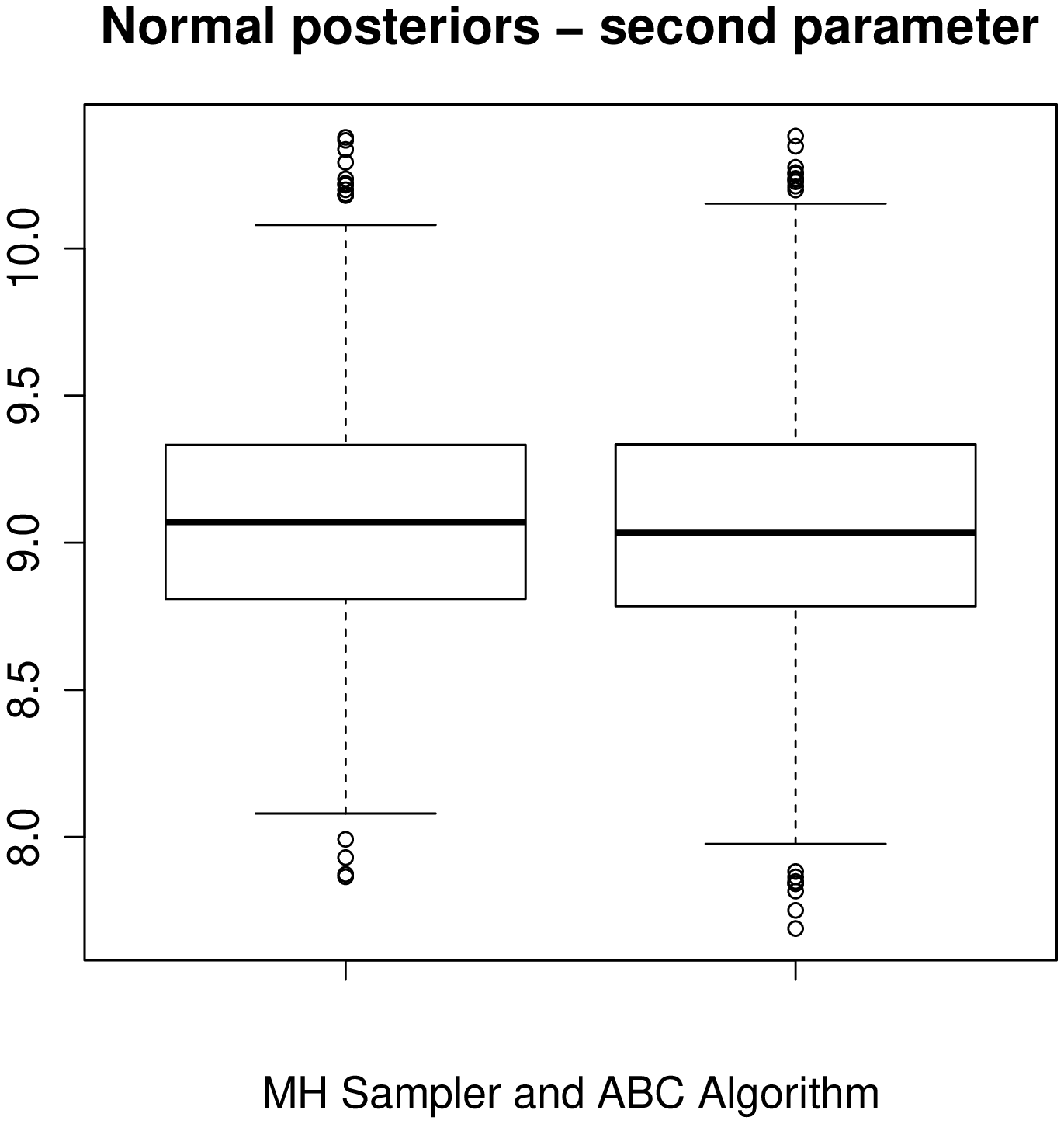} & \epsfxsize=6cm \epsffile{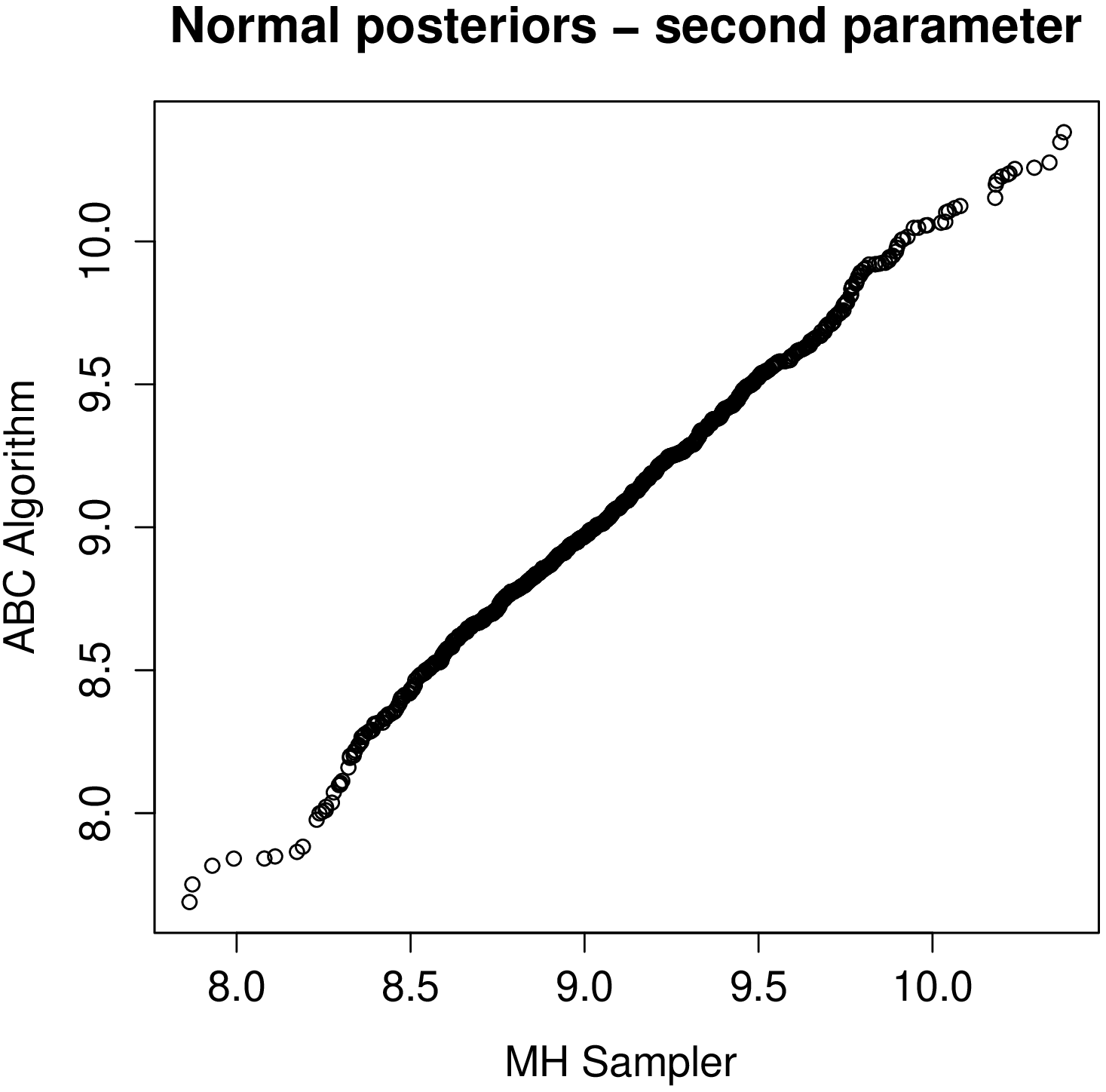}\\
\end{tabular}
\end{center}
\caption{Normal posterior sampling using a MH algorithm and the ABC Shadow procedure: boxplots and qqplots of the algorithms outputs.}
\label{resultsNormal}
\end{figure}

\noindent
Figure~\ref{resultsNormal} presents the results of the MH and the ABC Shadow algorithms used to sample from the Normal posterior~\eqref{posteriorNormal}. Table~\ref{tableNormal} shows the values of some summary statistics obtained using both methods. The posterior approximation looks satisfactory. It is important to notice, that the mean and the median posterior estimates are close to the maximum likelihood estimates of the model parameters.\\

\begin{table}[!htbp]
\begin{center}
\begin{tabular}{|c|c|c|c|c|c|c|} 
\hline
\multicolumn{7}{|c|}{Summary statistics for Normal posterior sampling}\\
\hline
Algorithm & $Q_{5}$ & $Q_{25}$ & $Q_{50}$ & $\bar{\theta}$ & $Q_{75}$ & $Q_{95}$ \\
\hline
MH $\theta_1$  & 1.60 & 1.69 & 1.75 & 1.76 & 1.82 & 1.92\\
ABC $\theta_1$ & 1.60 & 1.70 & 1.76 & 1.76 & 1.82 & 1.91\\
\hline
MH $\theta_2$  &  8.45 & 8.80 & 9.07 & 9.08 & 9.33 & 9.76\\
ABC $\theta_2$ &  8.35 & 8.78 & 9.03 & 9.06 & 9.33 &  9.83\\
\hline
\end{tabular}
\end{center}
\caption{Empirical quantiles and mean for the posterior of the Normal model.}
\label{tableNormal}
\end{table}

\begin{figure}[!htbp]
\begin{center}
\begin{tabular}{cc}
\epsfxsize=6cm \epsffile{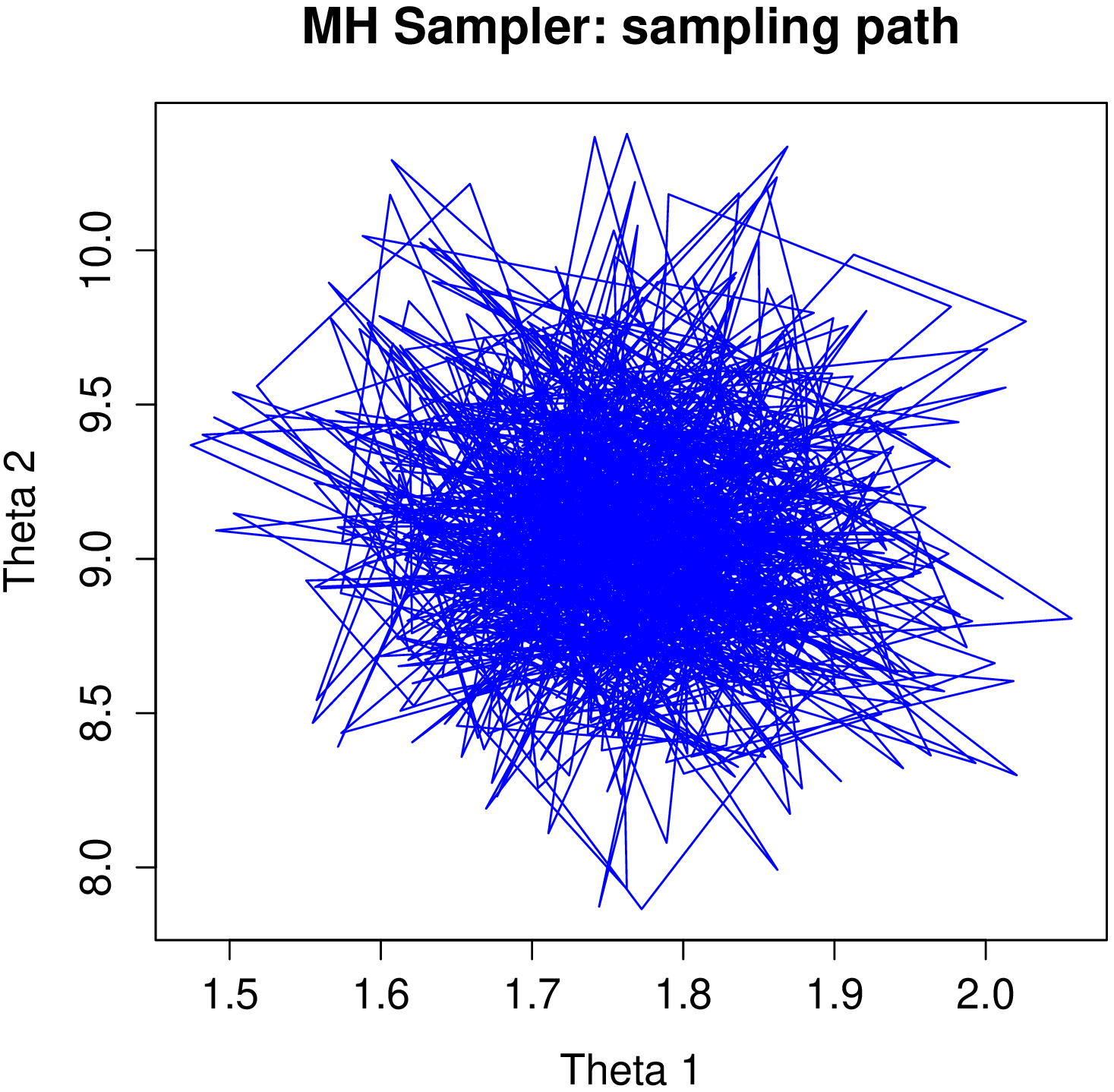} & \epsfxsize=6cm \epsffile{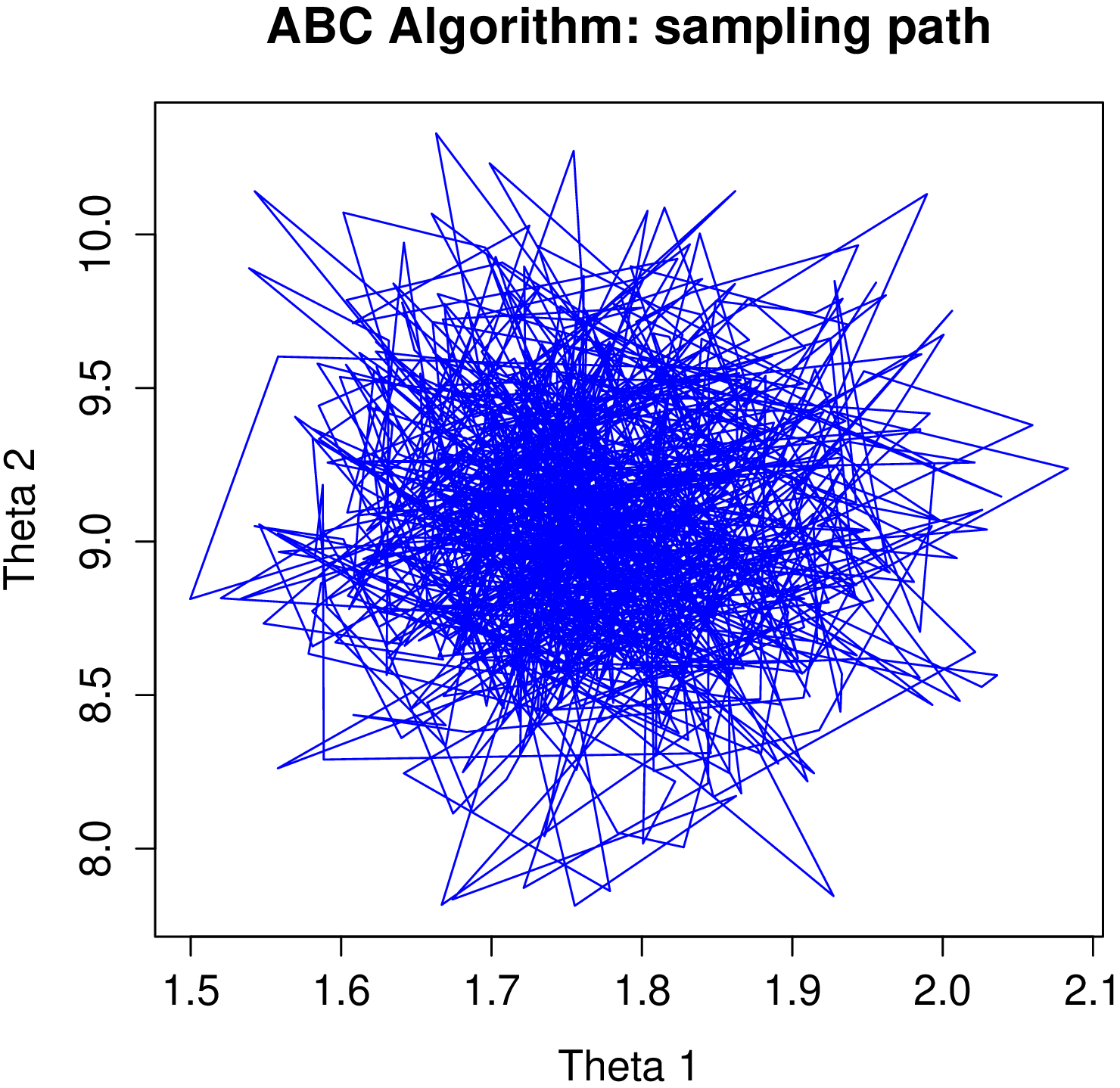} \\
\epsfxsize=6cm \epsffile{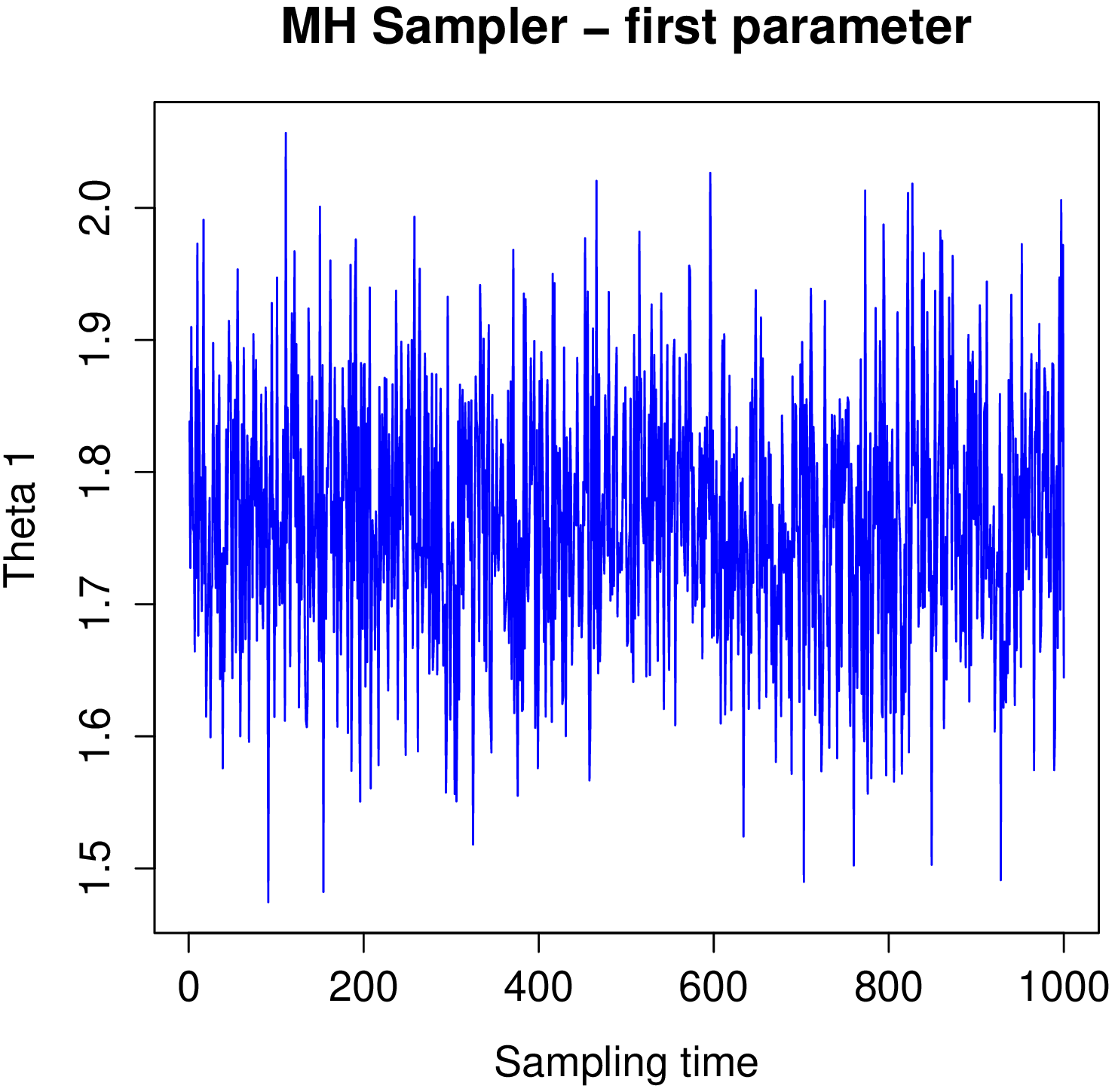} & \epsfxsize=6cm \epsffile{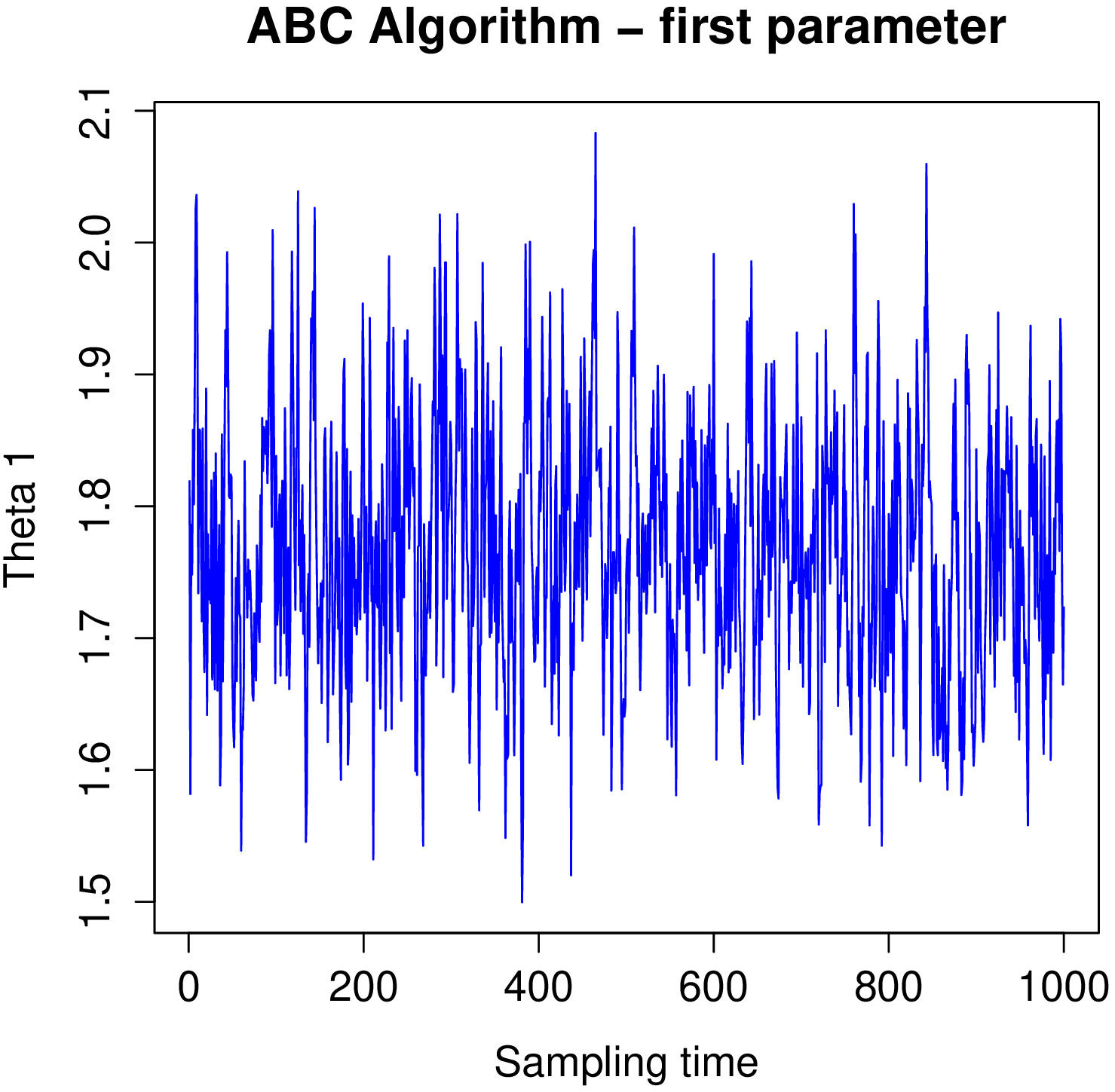}\\
\epsfxsize=6cm \epsffile{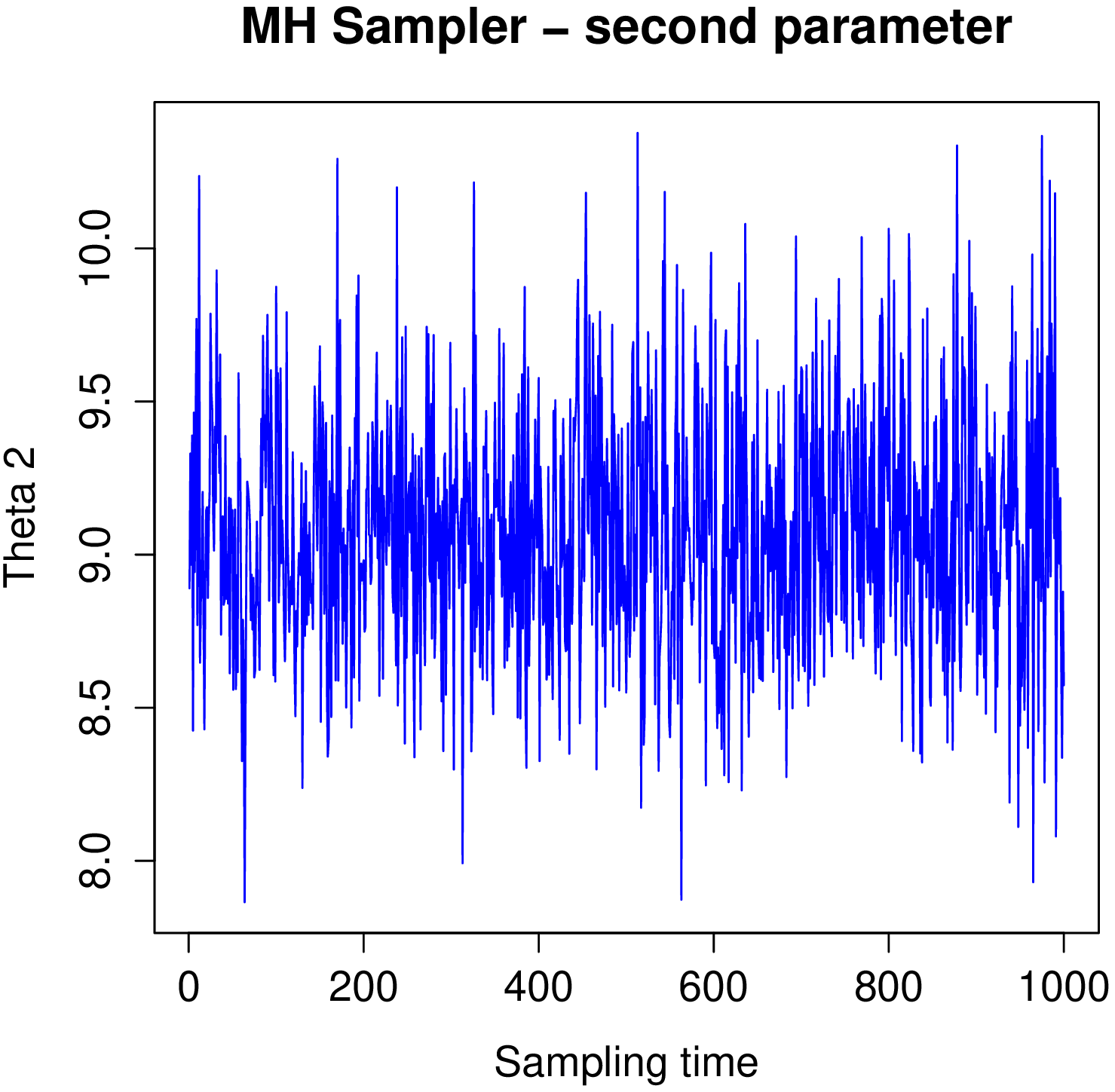} & \epsfxsize=6cm \epsffile{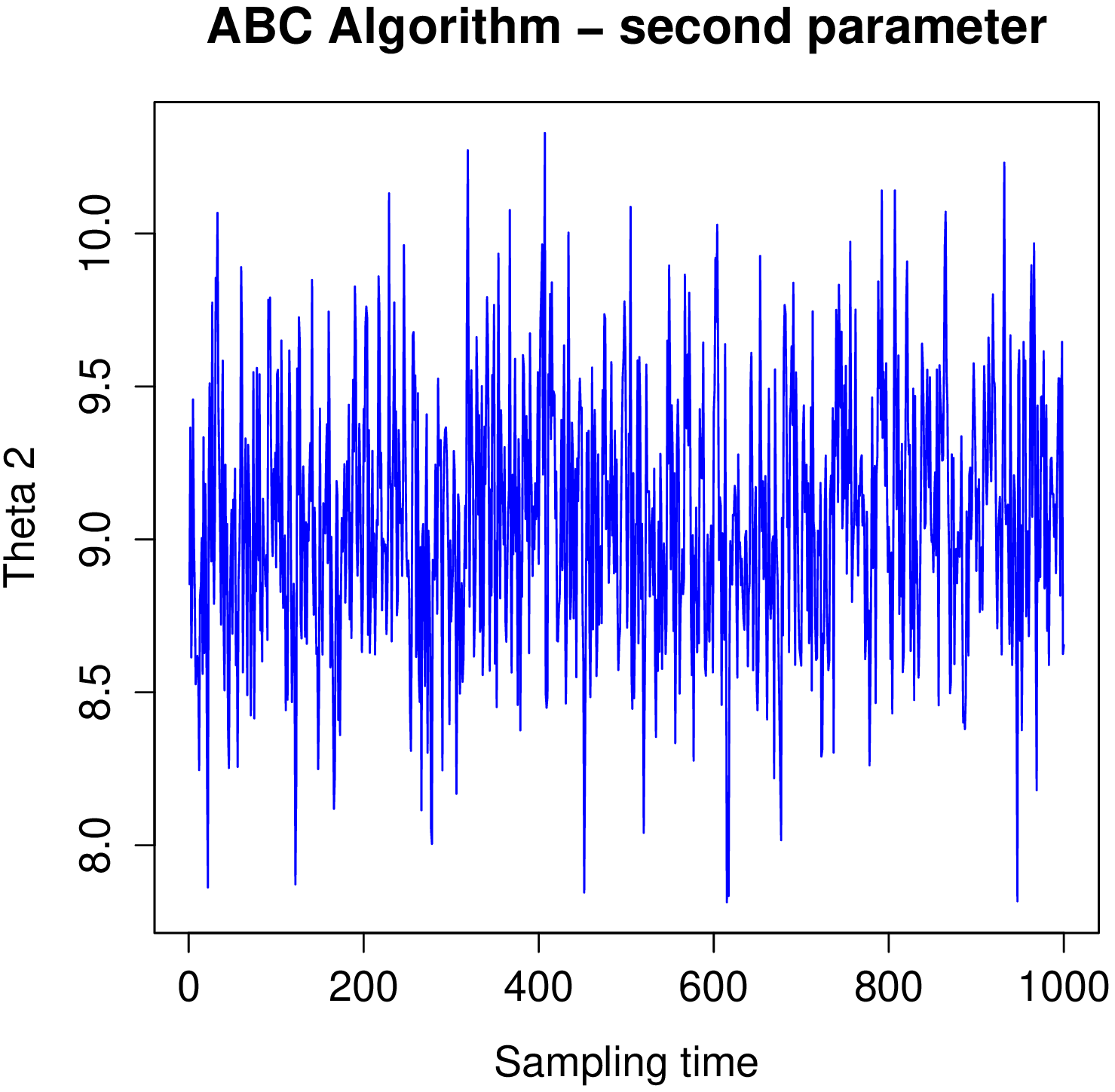}\\
\end{tabular}
\end{center}
\caption{Sample path for the Normal posterior. Left colum: the MH algorithm results - from top to bottom the joint parameter path, the $\theta_1$ time series and the $\theta_2$ time series. Right column: the ABC Shadow procedure - from top to bottom the joint parameter path, the $\theta_1$ time series and the $\theta_2$ time series.}
\label{pathNormal}
\end{figure}

\noindent
Figure~\ref{pathNormal} shows the sampling paths obtained using the MH and the ABC Shadow algorithms, respectively. In both cases, the sample paths exhibit a similar behaviour. The ABC chain indicates rather good mixing properties, in the sense that the chain does not get locked in a certain parameter space region and that it covers the parameter region given by the posterior distribution.\\

\noindent
Summary statistics obtained for the ABC algorithm run for different parameters $\Delta$ and initial conditions $\theta_0$ are presented in Table~\ref{analyseNormal}. For the rest of the parameters the Shadow algorithm was set as above. The obtained summary statistics show that if the values of $\Delta$ are too high, such as in $\Delta=(0.1,0.1)$, then the approximate posterior is rather far away from the true posterior. On the other way around, if $\Delta = (0.001,0.005)$ the summary statistics seem to indicate a good behaviour of the algorithm. Nevertheless, a visual investigation of the sample path confirmed as expected, that in this case, the simulated chain is slow mixing. Clearly, there is a compromise in choosing the $\Delta$ parameter. A too high value of $\Delta$ will make the algorithm diverge, since the distance between the approximate and the true posterior linearly depends on $\Delta$. An algorithm using a too small value of this parameter will produce samples that indeed may tend to approach slowly the true posterior distribution, but the mixing properties of the simulated chain are not guaranteed. The alternative intermediate values we have used for $\Delta$ show good values for the summary statistics, a proper mixing behaviour and independence of the initial conditions $\theta_0$.\\

\begin{table}[!htbp]
\begin{center}
\begin{tabular}{|c|c|c|c|c|c|c|} 
\hline
Parameters & \multicolumn{6}{|c|}{Summary statistics for Normal posterior sampling}\\
           & \multicolumn{6}{|c|}{$\theta_1$,$\theta_2$}\\
\hline
$\Delta$ , $\theta_0$ & $Q_{5}$ & $Q_{25}$ & $Q_{50}$ & $\bar{\theta}$ & $Q_{75}$ & $Q_{95}$ \\
\hline
$\Delta = (0.001,0.005)$ & 1.68 & 1.76 & 1.82 & 1.82 & 1.87 & 1.95\\
$\theta_0 = (2,9)$       & 8.41 & 8.81 & 9.13 & 9.09 & 9.38 & 9.83\\
\hline
$\Delta = (0.1,0.1)$ & -0.45 & 0.21 & 1.96 & 1.79 & 3.30 & 3.94\\
$\theta_0 = (2,9)$   &  5.87 & 6.57 & 7.03 & 7.06 & 7.52 & 8.31\\
\hline
$\Delta = (0.01,0.05)$ & 1.58 & 1.69 & 1.76 & 1.76 & 1.83 & 1.92\\
$\theta_0 = (2,9)$     & 8.42 & 8.77 & 9.03 & 9.07 & 9.33 & 9.82\\
\hline
$\Delta = (0.005,0.025)$ & 1.61 & 1.70 & 1.76 & 1.76 & 1.83 & 1.93\\
$\theta_0 = (2,9)$       & 8.44 & 8.78 & 9.05 & 9.07 & 9.32 &  9.74\\
\hline
$\Delta = (0.005,0.025)$ & 1.60 & 1.70 & 1.76 & 1.76 & 1.82 & 1.91\\ 
$\theta_0 = (10,20)$     & 8.35 & 8.78 & 9.03 & 9.06 & 9.33 & 9.83\\
\hline
$\Delta = (0.005,0.025)$ & 1.61 & 1.70 & 1.77 & 1.76 & 1.83 & 1.92\\
$\theta_0 = (-10,1)$     & 8.39 & 8.80 & 9.06 & 9.08 & 9.36 & 9.80\\
\hline
\end{tabular}
\end{center}
\caption{ABC Shadow behaviour depending on the $\Delta$ parameter and $\theta_0$ initial conditions.}
\label{analyseNormal}
\end{table}

\subsection{Posterior approximation for the analysis of spatial patterns}
\noindent
The Shadow algorithm is applied here to the statistical analysis of three spatial data sets. The first two data sets consist of several patterns which are realisations of a Strauss model and a Candy model, respectively. The third one is a a single point pattern which is real cosmological data~\cite{StoiEtAl15}.\\

\subsubsection{Strauss model patterns}
The Strauss model~\cite{KellRipl76,Stra75} simulates random patterns made of repulsive points. The probability density of the process is
\begin{equation}
p(\yy|\theta) \propto \beta^{n(\yy)}\gamma^{s_{r}(\yy)} = \exp \langle n(\yy) \log \beta  + s_r(\yy)\log \gamma \rangle.
\label{straussModel}
\end{equation}
Here $\yy$ is a point pattern in the window $W$, while $t(\yy) = (n(\yy),s_{r}(\yy))$ and $\theta = (\log \beta, \log \gamma)$ are the sufficient statistics and the model parameters vectors, respectively. The sufficient statistics components $n(\yy)$ and $s_{r}(\yy)$ represent each, the number of points in $W$ and the number of pairs of points at a distance closer than $r$. The $\beta$ parameter controls the number of points in a configuration, while the $\gamma$ parameter controls the relative position of the points in a configuration. The model is well defined for $\beta > 0$ and $\gamma \in ]0,1]$. If $\gamma = 1$ then the point process is a Poisson point process, since no interaction between points exists. If $\gamma \in (0,1)$ the points in a configuration exhibit a repulsion interaction. The range parameter $r$ is usually considered known. It cannot be taken directly into account in our framework since the likelihood $p(\yy|\theta)$ is not a differentiable function of it.\\

\noindent
Here, it is not possible to make a direct comparison that certifies the quality of the results given by the ABC Shadow algorithm. Nevertheless, for some particular situations, it may be investigated whether the inference from the approximated posterior distribution is close to the inference performed using classical methods.\\

\noindent
Let us assume that a random sample of size $m$ made of independent realisations of a Strauss model with parameter $\theta$ is observed. It is easy to check following~\cite{Lies00,MollWaag04}, that the maximum likelihood estimate $\widehat{\theta}$ satisfies the equation
\begin{equation*}
\sum_{i=1}^{m}t(\yy_i) - m \EE_{\widehat{\theta}}t(\XX) = 0,
\end{equation*}
with $\yy_1,\ldots,\yy_m$ the $m$ point patterns forming the sample. The sum and the expectation are computed componentwise. Clearly, the maximum likelihood estimate approaches the true model parameters whenever the sample size $m$ increases. Since the considered parameter domain is a compact set $\Theta \subset \RR^{r}$, the likelihood based inference and the one based on a posterior distribution build with uniform priors on $\Theta$, are equivalent.\\

\noindent
This property can be used to test the ABC Shadow algorithm: within this context, the maximum of the approximate posterior should be close to the maximum likelihood estimate and the true model parameters. For saving computational time, instead of observing $m$ independent samples, an average pattern with sufficient statistics given by $\frac{1}{m}\sum_{i=1}^{m}t(\yy_i)$ can be considered.\\

\noindent
In the following, the previous ideas are used to test our method. For this purpose, the Strauss model on the unit square $W=[0,1]^2$ and with density parameters $\beta=100$, $\gamma=0.2$ and $r=0.1$, was considered. This gives for the parameter vector of the exponential model $\theta=(4.60,-1.60)$. The exact algorithm CFTP~(see Chapter 11 in \cite{MollWaag04}) was used to get $1000$ samples from the model and to compute the empirical means of the sufficient statistics $\bar{t}(\yy) = (\bar{n}(\yy),\bar{s_{r}}(\yy))=(34.33,5.31)$. The ABC Shadow algorithm was run using this sufficient statistics vector as observed data. The prior density $p(\theta)$ was the uniform distribution on the interval $[3.5,5.5] \times [-5,0]$. Each time, the auxiliary variable was sampled using $100$ steps of a MH dynamics~\cite{Lies00,MollWaag04}. The $\Delta$ and $n$ parameters were set to $(0.01,0.01)$ and $(200,200)$. The algorithm was run for $10^6$ iterations. Samples were kept every $10^3$ steps. This gave a total of $1000$ samples.\\

\begin{figure}[!htbp]
\begin{center}
\begin{tabular}{cc}
\epsfxsize=6cm \epsffile{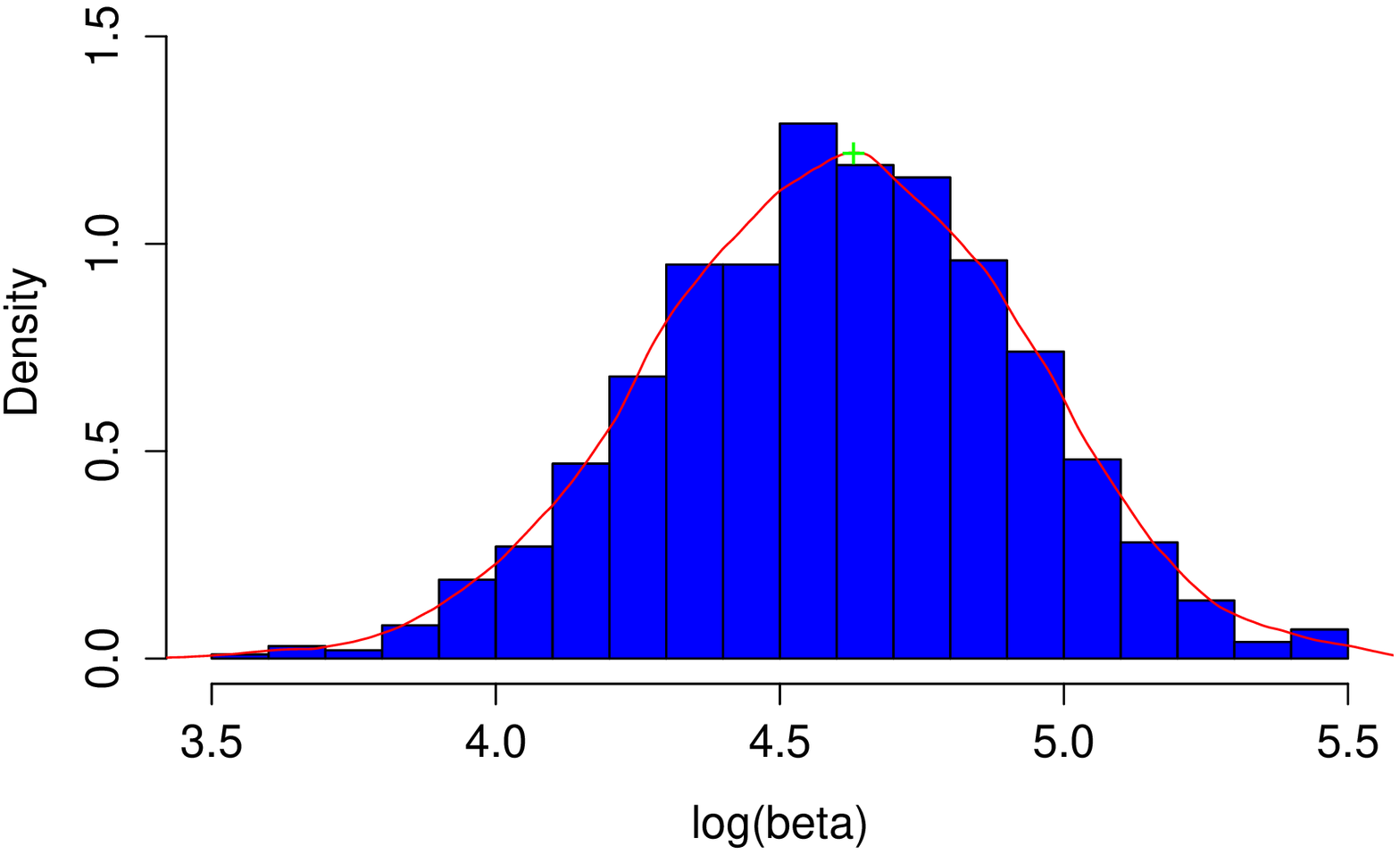} &
\epsfxsize=6cm \epsffile{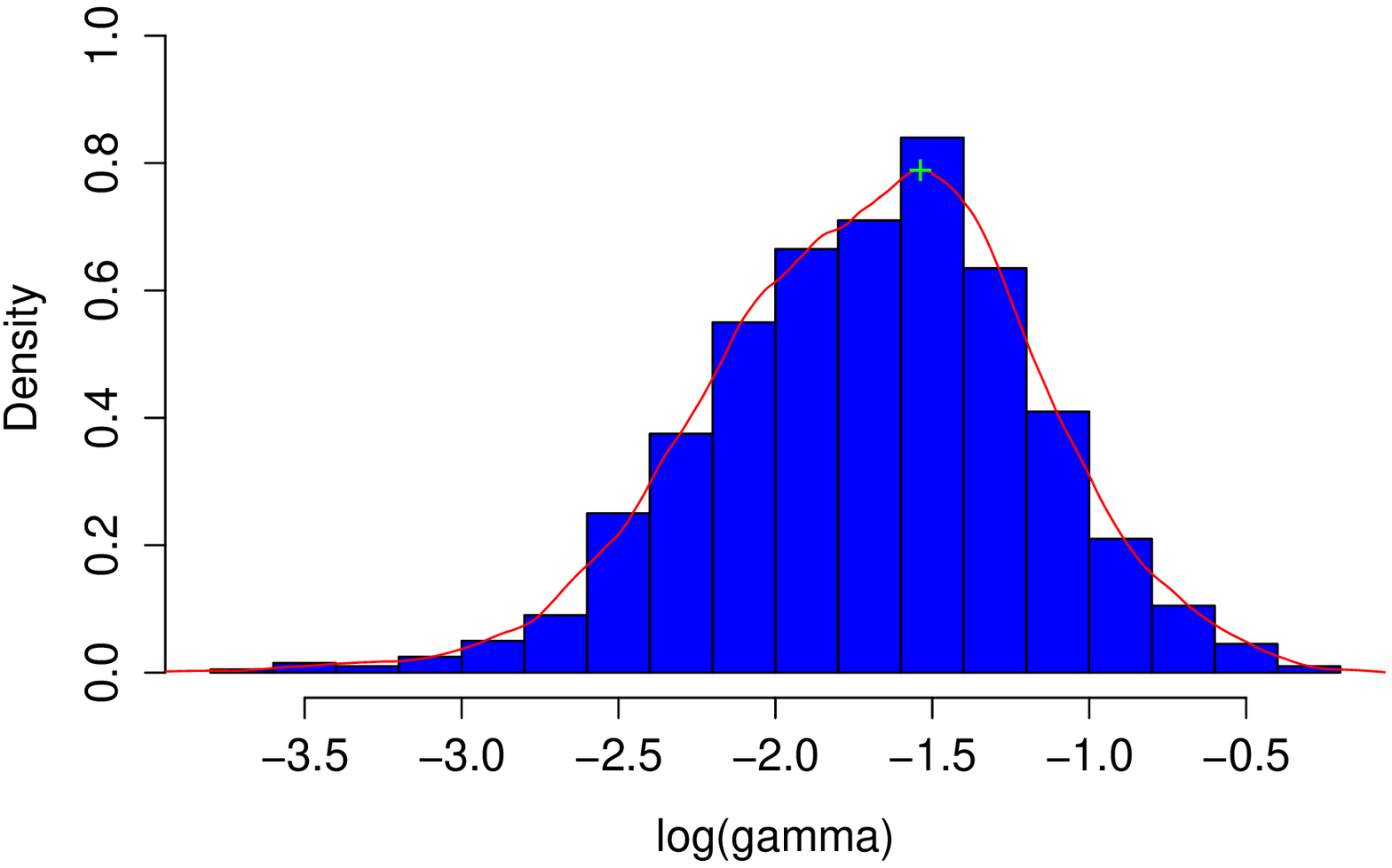}
\end{tabular}
\end{center}
\caption{Strauss posterior sampling using the ABC procedure: histograms and kernel density estimates of the marginals.}
\label{resultsStrauss}
\end{figure}

\noindent
Figure~\ref{resultsStrauss} shows the histograms obtained using the ABC algorithm for sampling the posterior density given by the Strauss model~\eqref{straussModel}. The kernel density estimates are superimposed on these histograms. The MAP (Maximum a posteriori) is computed by taking the maximum of these estimated densities. The obtained value is $\widehat{\theta} = (4.63,-1.53)$. The obtained median and mean posterior estimates were $(4.606,-1.669)$ and $(4.603,-1.700)$, respectively. All these values are close to the true values of the parameters. Note that in this case, the true maximum likelihood estimates are also equal to the true model parameters. Nevertheless, nothing can be said concerning the general shape of the distribution.\\

\subsubsection{Candy model patterns}
The Candy model is a marked point process that simulates connected segments~\cite{LiesStoi03}. The model was applied with success in image analysis and cosmology~\cite{StoiDescZeru04,StoiEtAl15}.\\

\noindent
Here, a segment or a marked point $y=(w,\xi,l)$ is given by its centre $w$, its orientation $\xi$ and its length $l$. The orientation mark is a uniform random variable on $M=[0,\pi)$, while the length is a fixed value. The Candy model we consider is given by the following probability density
\begin{equation}
p(\yy|\theta) \propto \exp \langle \theta_{d} n_{d}(\yy) + \theta_{s} n_s(\yy) + \theta_{f} n_f(\yy) + \theta_{r} n_{r}(\yy) \rangle
\label{candyModel}
\end{equation}
with $\theta = (\theta_d,\theta_s,\theta_f,\theta_r)$ and $t(\yy) = (n_{d}(\yy), n_{s}(\yy), n_{f}(\yy), n_{r}(\yy))$ the parameters and sufficient statistics vectors, respectively. The parameter $\theta_d$ controls the number $n_d$ of segments connected at both of its extremities or doubly connected, the parameter $\theta_s$ controls the number $n_s$ of segments connected at only one of its extremities or singly connected and the parameter $\theta_f$ controls the number $n_f$ of segments that are not connected or free. The connectivity interaction of two segments is based on the relative position of their extremities and also on their relative orientation. Two segments with only one  pair of extremities situated within connection range $r_c$ and with absolute orientation difference lower than a curvature parameter $\tau_c$ are connected. The parameter $\theta_r$ controls the number $n_r$ of pairs of segments that are too close and not orthogonal. This interaction controls the segments tendency to form clusters or to overlap. The orthogonality is controlled through a curvature parameter $\tau_r$. For more details and a complete description of the Candy model interactions and properties we recommend~\cite{LiesStoi03}.\\

\begin{figure}[!htbp]
\begin{center}
\epsfxsize=10cm \epsffile{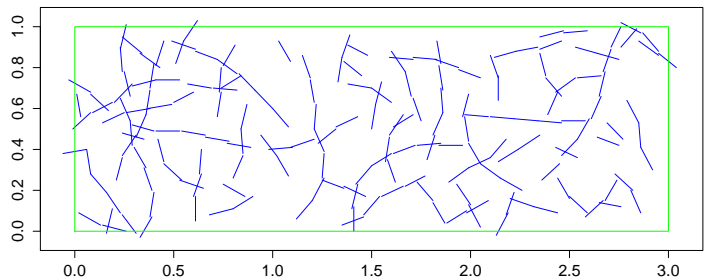} 
\end{center}
\caption{Realisation of the Candy model.}
\label{sampleCandy}
\end{figure}

\noindent
Figure~\ref{sampleCandy} shows a realisation of the Candy model on $W=[0,3] \times [0,1]$. The segment length is $l=0.12$, the connection range is $r_c = 0.01$, and the curvature parameters are $\tau_c=\tau_r=0.5$ radians. The model parameters are $\theta_d = 10$, $\theta_s=7$, $\theta_f=3$ and $\theta_r = -1$. It can be observed that with these parameters a rather connected pattern of segments is formed.\\

\noindent
The same strategy as for the Strauss model was used to test the ABC Shadow algorithm. An adapted MH algorithm~\cite{LiesStoi03} was used to obtain $1000$ samples of the previous model and to compute the vector of the empirical means of the sufficient statistics $\bar{t}(\yy) = (\bar{n}_{d}(\yy) = 51.10, \bar{n}_{s}(\yy) = 101.06, \bar{n}_{f}(\yy) = 19.97, \bar{n}_{r}(\yy) = 72.89)$. These statistics were considered as the observed data for our ABC Shadow algorithm. The prior density $p(\theta)$ was the uniform distribution on the interval $[2,12] \times [2,12] \times [2,12] \times [-7,0]$. Each time, the auxiliary variable was sampled using $2000$ steps of an adapted MH dynamics. The $\Delta$ and $n$ parameters were set to $(0.01,0.01,0.01,0.01)$ and $(500,500,500,500)$. The algorithm was run for $10^6$ iterations. Samples were kept every $10^3$ steps. This gave a total of $1000$ samples.\\

\begin{figure}[!htbp]
\begin{center}
\begin{tabular}{cc}
\epsfxsize=6cm \epsffile{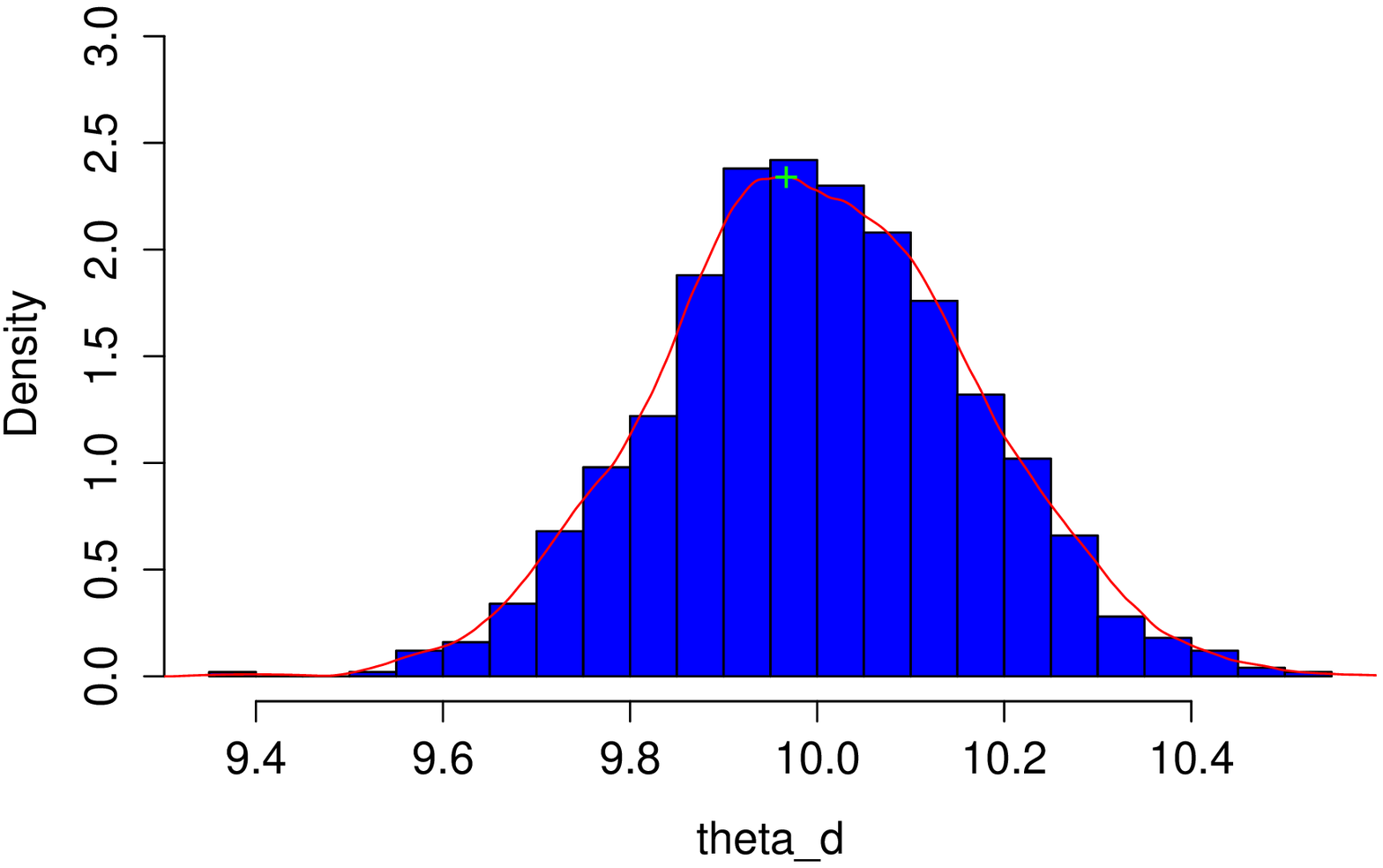} &
\epsfxsize=6cm \epsffile{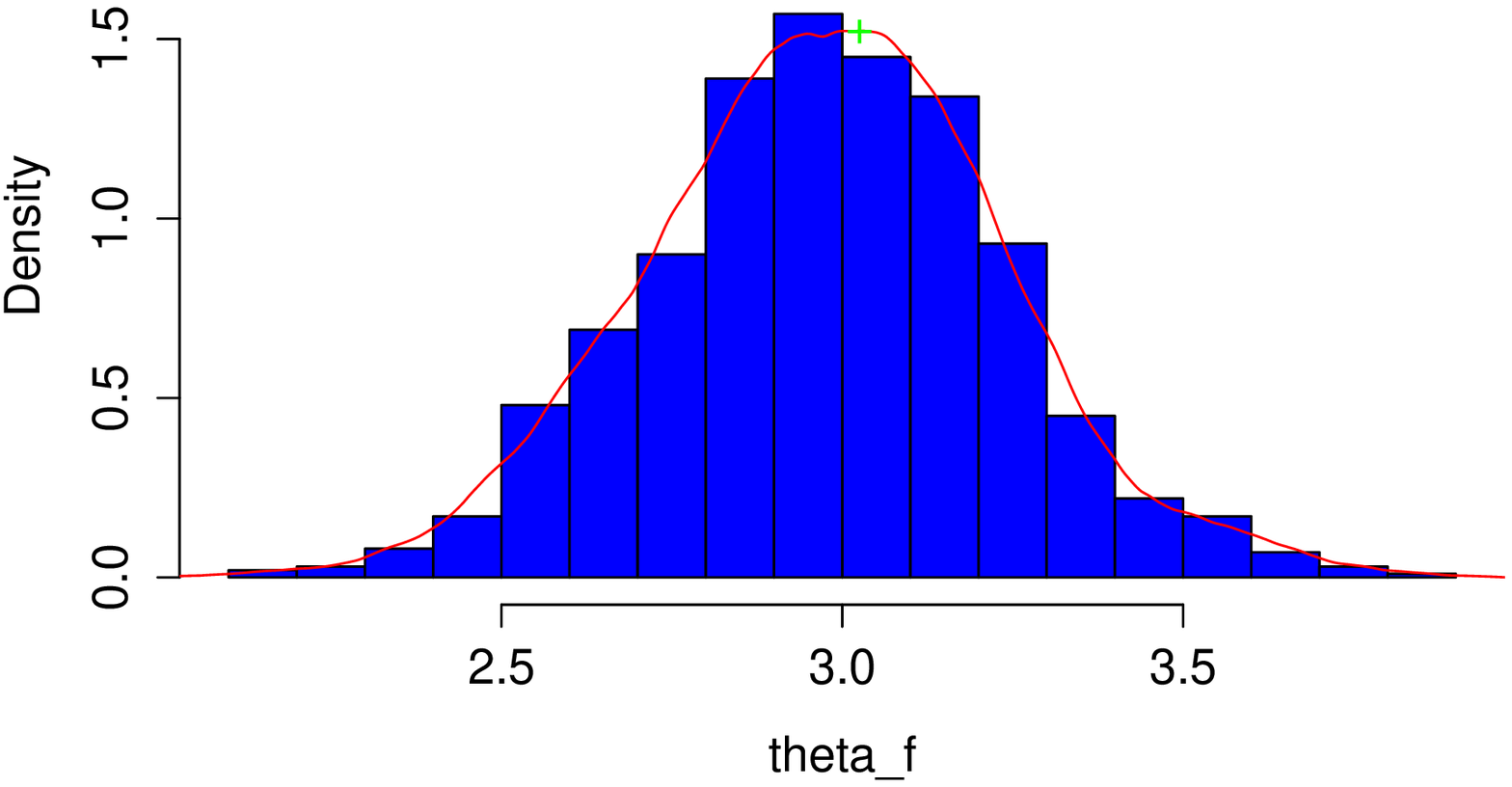}\\
\epsfxsize=6cm \epsffile{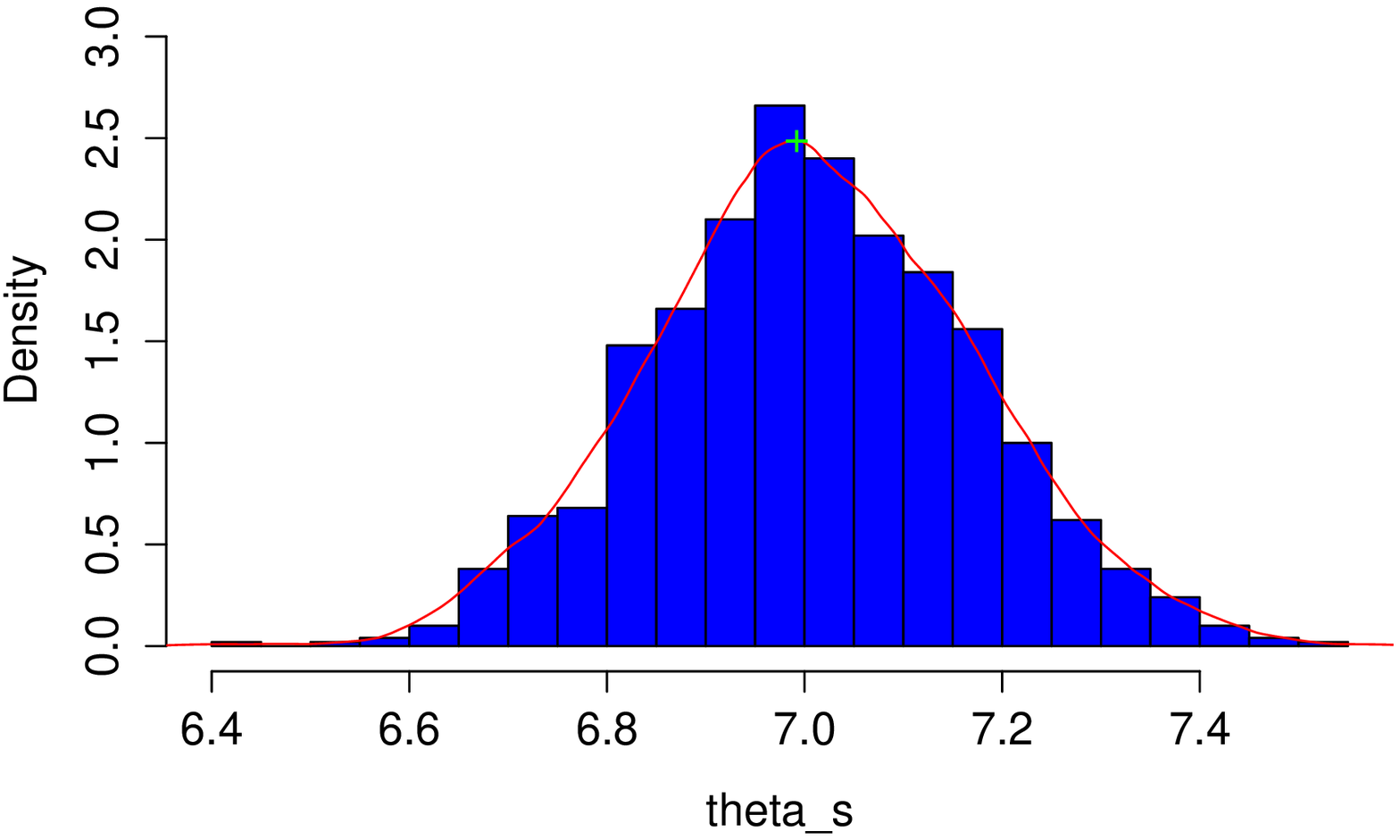} &
\epsfxsize=6cm \epsffile{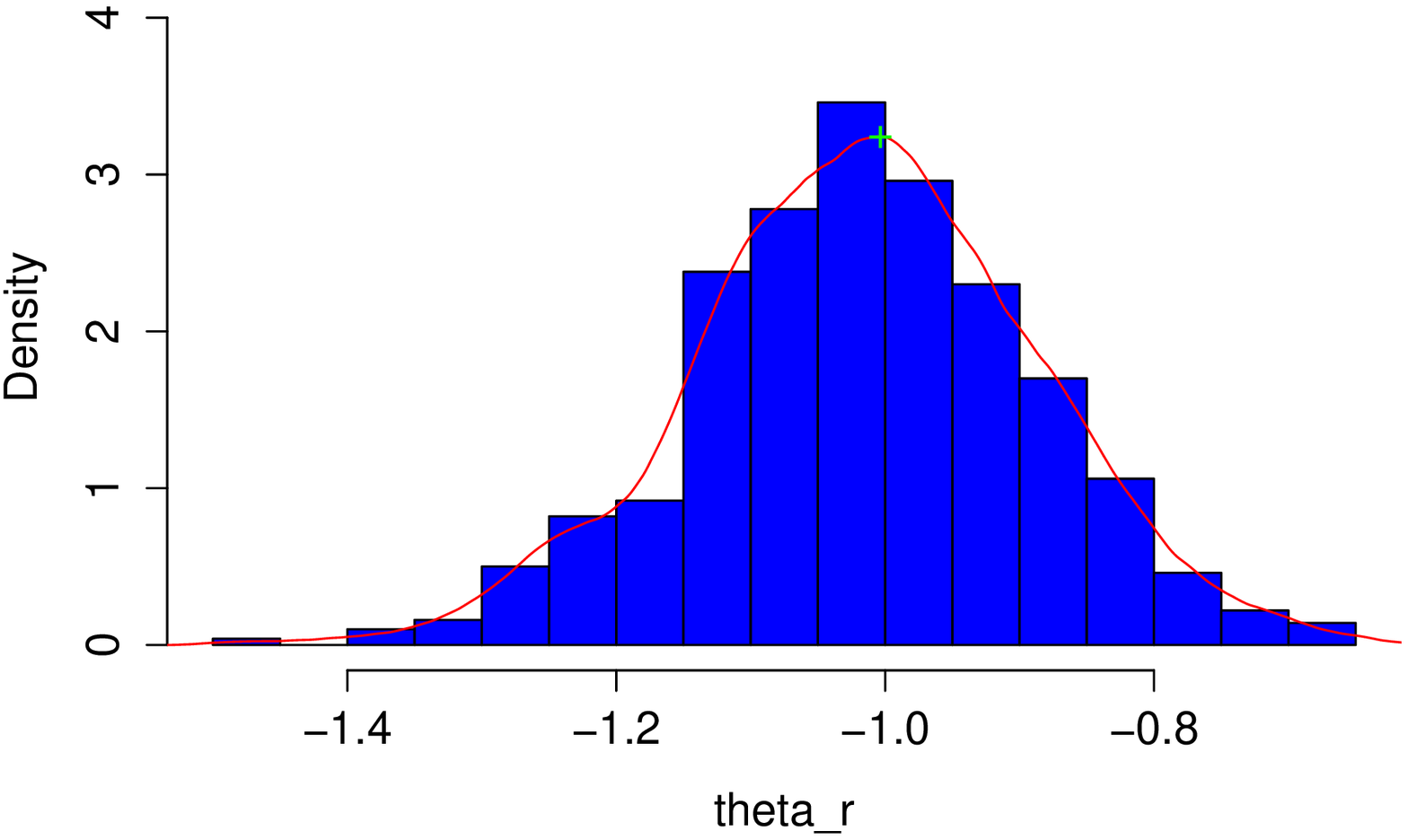}\\
\end{tabular}
\end{center}
\caption{Candy posterior sampling using the ABC procedure: histograms and kernel density estimates of the marginals.}
\label{resultsCandy}
\end{figure}

\noindent
Figure~\ref{resultsCandy} shows the histograms obtained using the ABC algorithm for sampling the posterior density given by the Candy model~\eqref{candyModel}. The kernel densities estimates are superimposed on these histograms. The MAP is computed by taking the maximum of these estimated densities. The obtained value is  $\widehat{\theta}=(9.96, 6.99, 3.02,-1.00)$. The obtained median and mean estimates were $(9.995, 7.005, 2.977,-1.014)$ and $(9.998, 7.008, 2.975,-1.018)$, respectively. Again, the obtained values are close to the true parameter values.\\
 
\subsubsection{Real data application}
\begin{figure}[!htbp]
\begin{center}
\epsfxsize=10cm \epsffile{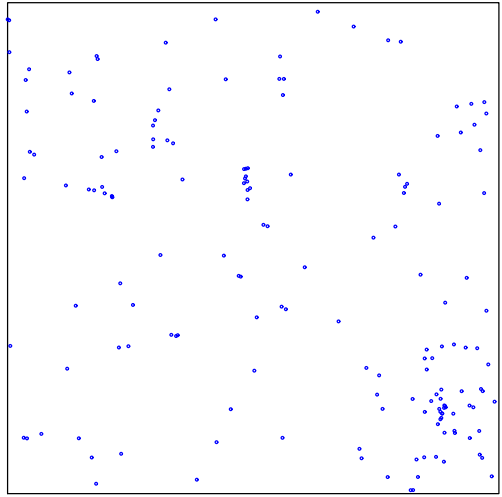} 
\end{center}
\caption{Sample of a cosmological data set. The points represent galaxy positions in a region of our Universe.}
\label{galaxyData}
\end{figure}

\noindent
The galaxies in our Universe are not uniformly distributed. They tend to form structures such as filaments, sheets and clusters~\cite{MartSaar02}. Figure~\ref{galaxyData} shows a sub-region from a two-dimensional cosmological catalogue~\cite{StoiEtAl15}. The data set contains a number of $163$ points that represent the galaxy positions in the considered region. The observation domain is scaled to a square window $W = [0,1] \times [0,1]$. The visual inspection of the data set indicates a clustering tendency of the pattern. In the following, it will be shown how the ABC Shadow algorithm can be used to test this hypothesis.\\

\noindent
Such point patterns are typically analysed using summary statistics~\cite{Lies00,MollWaag04}. For a stationary point process $\YY$ of intensity $\lambda >0$, some of the most used summary statistics are:
\begin{itemize}
\item the $K$ function, which is proportional to the expectation under the reduced Palm distribution of $N(b(o,u))$, the number of points in a ball of radius $u$ centred at the origin $o$, that is a point of the considered stationary process which is {\it not} counted
\begin{equation*}
\lambda K(u) = \EE_{o}^{!}\left[ N(b(o,u))\right],
\end{equation*}
\item the nearest neighbour distance distribution function
\begin{equation*}
G(u) = \PP_{o}^{!}(N(b(o,u) > 0))
\end{equation*}
with $\PP_{o}^{!}$ the reduced Palm distribution,
\item the empty space function
\begin{equation*}
F(u) = \PP(\YY \cap b(o,u) \neq \emptyset) 
\end{equation*}
with $\PP$ the distribution of $\YY$,
\item the $J$ function which compares nearest neighbour to empty distances
\begin{equation*}
J(u) = \frac{1-G(u)}{1-F(u)}
\end{equation*}
for all $u > 0$ such that $F(u) < 1$.
\end{itemize}
These summary statistics have exact formulas for the stationary Poisson point process with intensity $\lambda$~:
\begin{eqnarray*}
K(u) & = & \pi u^2, \nonumber \\
F(u) & = &  1 - \exp[-\lambda \pi u^2], \nonumber\\
G(u) & = & F(u), \nonumber \\
J(u) & = & 1.
\end{eqnarray*}
So, for a given point pattern, the estimators of these statistics  give indications about how far from the Poisson point process, the considered pattern is. For the $K$ and $G$ statistics, higher values of the observed summaries than the theoretical Poisson statistics suggest clustering of the pattern, while the lower ones suggest repulsion. For the $F$ and $J$ statistics, it is the other way around: lower observed values than the theoretical Poisson ones indicate clustering, while the higher values recommend repulsion. To better answer this type of questions for an observed point pattern, an envelope test can be done. This test simulates a Poisson point process with intensity estimated from the observed pattern, and for each simulated pattern the summary statistics are computed. At the end of the simulation, Monte Carlo confidence envelopes are created. So, if the observed summary statistics belong to the envelope, the null hypothesis that the observed point pattern is the realization of a Poisson process is not rejected, hence the point pattern exhibits a completely random structure. For a thorough presentation of the point process summary statistics and of the envelope tests based on them, we recommend~\cite{Lies00,MollWaag04} and the references therein.\\

\begin{figure}[!htbp]
\begin{center}
\begin{tabular}{cc}
\epsfxsize=6cm \epsffile{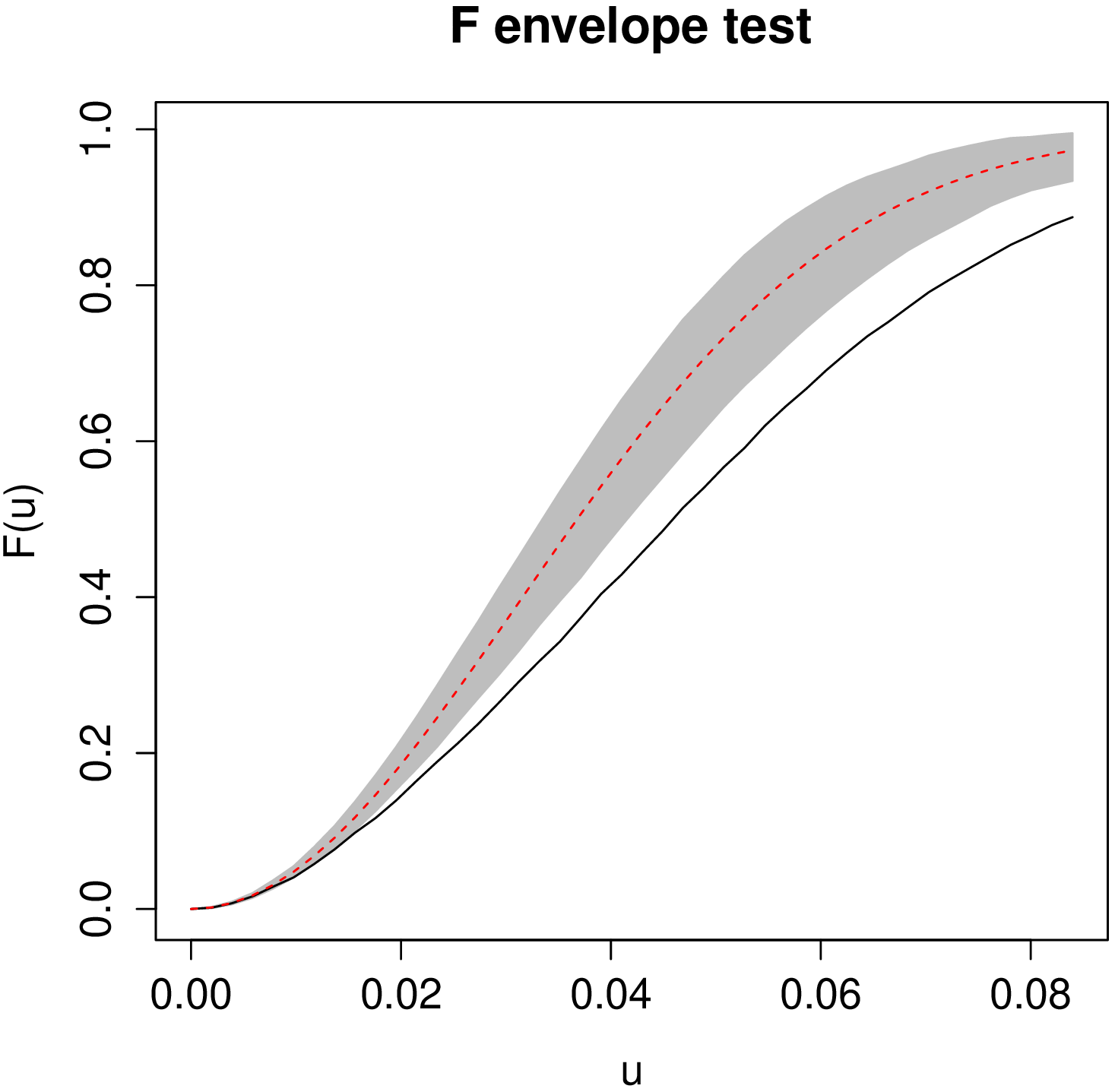} &
\epsfxsize=6cm \epsffile{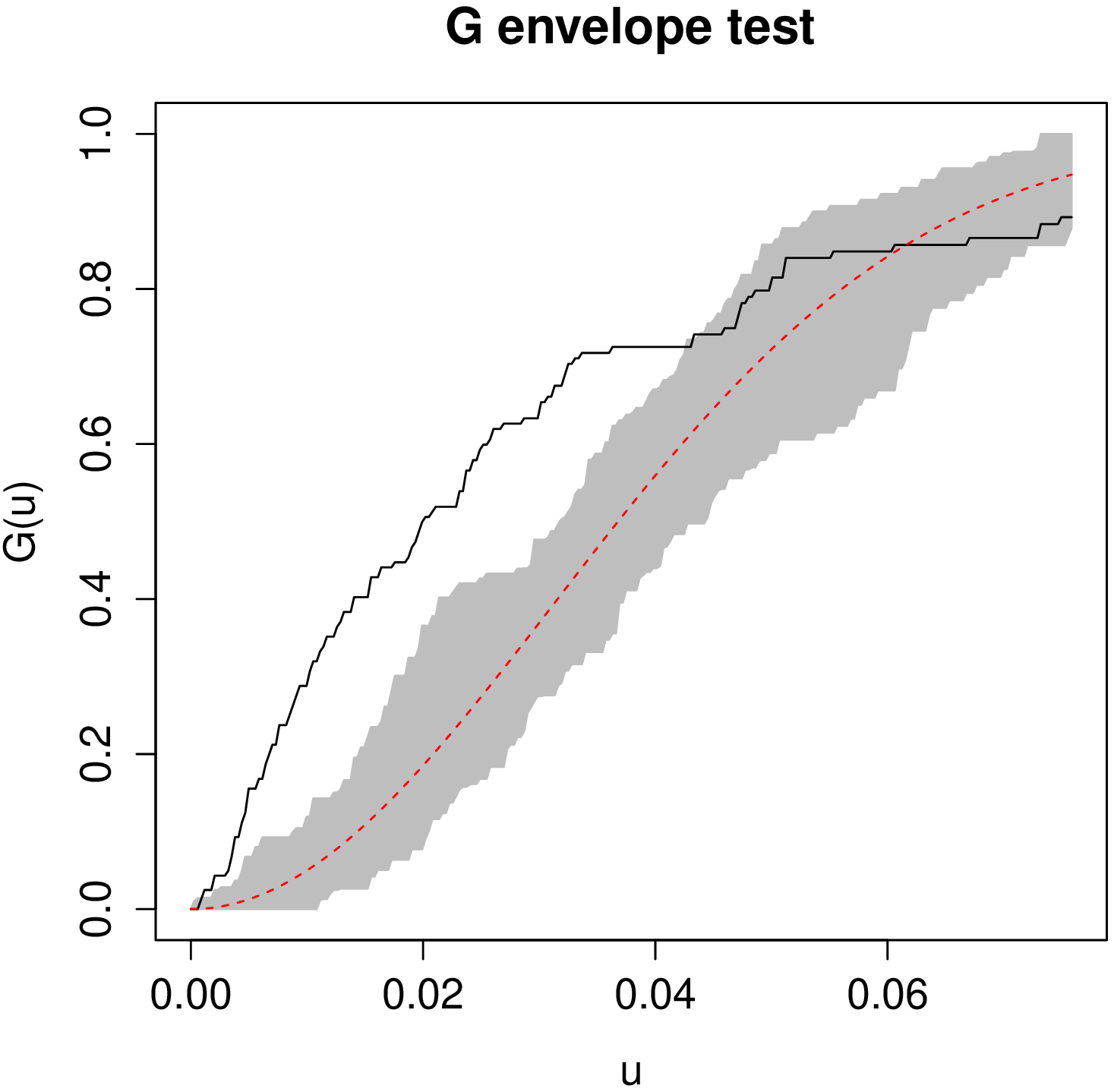}\\
\epsfxsize=6cm \epsffile{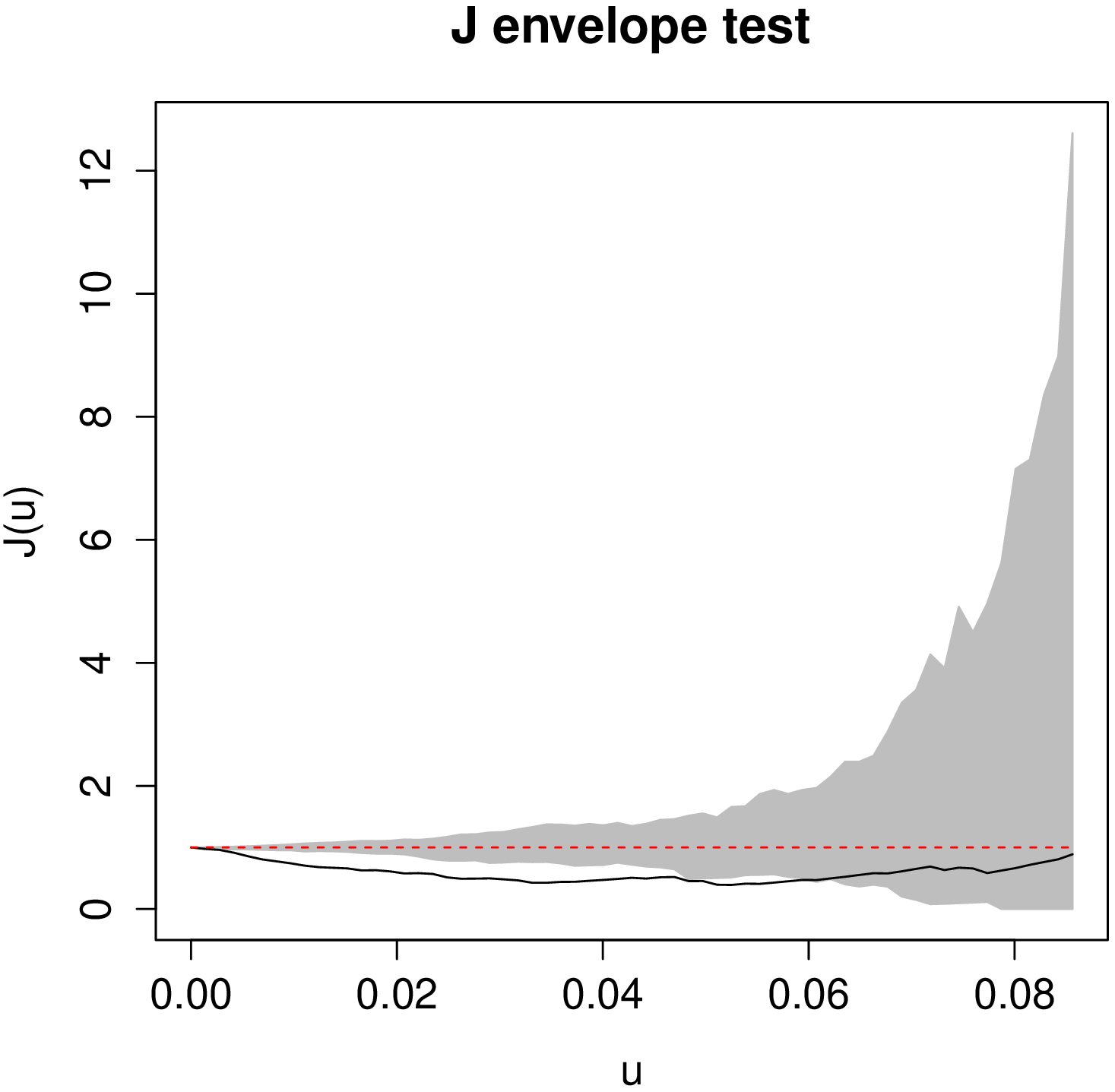} &
\epsfxsize=6cm \epsffile{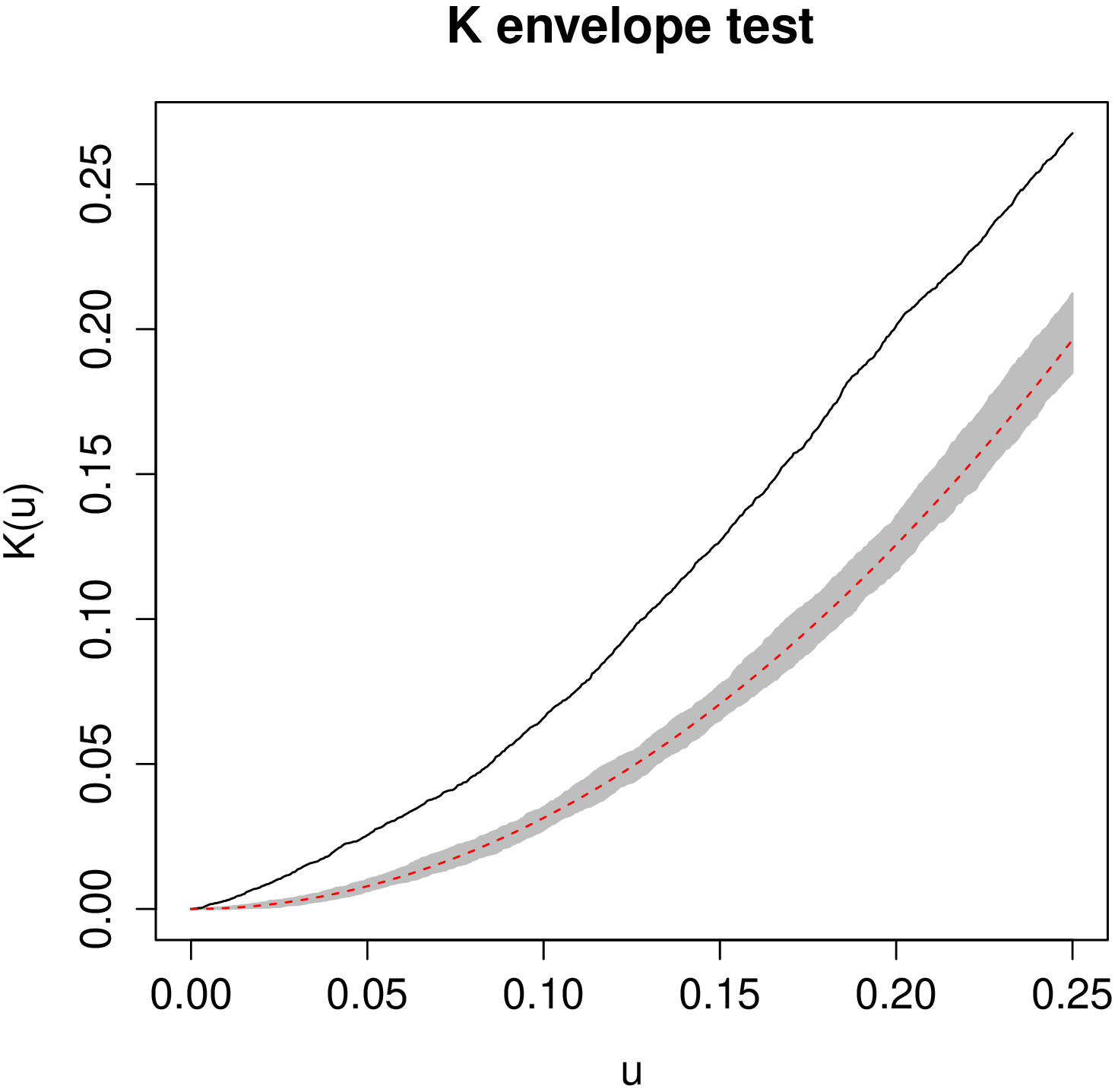}\\
\end{tabular}
\end{center}
\caption{Summary statistics envelope test ($100$ simulations) for the galaxy pattern: shaded (gray) region - the Monte Carlo envelopes, dotted line (red) - the theoretical statistics, continuous (black) line - the observed statistics. These results (estimators computations and simulations) were obtained using the {\bf R spatstat} library~\cite{BaddTurn05}.}
\label{resultsEnvelope}
\end{figure}

\noindent
Figure~\ref{resultsEnvelope} shows the result of the envelope for the galaxy point pattern. All the statistics exhibit the clustered character of the pattern. Despite the fact this very useful analysis result may be considered reliable, there is no information concerning the ``strength'' of the clustering character. This may be achieved by assuming the pattern as the realisation of known model and estimate its parameters. Still, parametric inference using these summary statistics should be done with care, since it may be misleading. For instance, the authors in~\cite{BaddSilv84,BedfBerg97} proved the existence of point processes that are not Poisson, but exhibiting $K(u) = \pi u^2$ or $J(u) = 1$. Another possibility to be considered is to perform pseudo-likelihood (PL) or Monte Carlo maximum likelihood (MCML) estimation. Now, the PL estimation lacks of precision if the objects interaction is too strong, while the MCML computation requires re-sampling of the model whenever the initial conditions are too far away from the true maximum likelihood. As an alternative to these, we use posterior sampling using the ABC Shadow algorithm.\\

\noindent
The area-interaction process was introduced by~\cite{BaddLies95} in order to produce cluster point patterns.  The probability density of the process is
\begin{equation}
p(\yy|\theta) \propto \beta^{n(\yy)}\gamma^{a_{r}(\yy)} = \exp \langle n(\yy) \log \beta  + a_r(\yy)\log \gamma \rangle.
\label{straussModel}
\end{equation}
Here $\yy = \{w_1,w_2,\ldots,w_{\kappa}\}$ is a point pattern in the finite window $W$, while $t(\yy) = (n(\yy),a_{r}(\yy))$ and $\theta = (\log \beta, \log \gamma)$ are the sufficient statistics and the model parameters vectors, respectively. The statistic $n(\yy)$ represents the number of points in a configuration. The statistic $a_{r}(\yy)$ is given by
\begin{equation*}
a_{r}(\yy) = -\frac{\nu\left[A_{r}(\yy))\right]}{\pi r^2} = - \frac{\nu\left[\cup_{i=1}^{n}b(w_i,r)\right]}{\pi r^2}
\end{equation*}
and it is proportional to the surface of the union disks of radius $r$ centred in the points of $\yy$. It can be interpreted also as the number of disks needed to ``cover'' the area of the disk pattern $A_{r}(\yy) = \cup_{i=1}^{n}b(w_i,r)$. Here, as outlined by~\cite{BaddLies95}, there is a direct link between the statistic $a_{r}(\yy)$ and the "reduced sample" estimate of the empty space function $\widehat{F}(r)$, in the sense that the vector $(n(\yy),\widehat{F}(r))$ is also a sufficient statistics vector. As for the Strauss model, the $\beta$ parameter controls the number of points in a configuration, while the $\gamma$ parameter controls the relative position of the points in a configuration. The model is well defined for $\beta,\gamma > 0$. If $\gamma = 1$ then the point process is a Poisson point process, and no interaction between points exists. If $\gamma > 1$ then the configurations inducing disk patterns with large areas will be penalised, hence producing a clustering effect of the points in the pattern. If $\gamma < 1$, then the configurations inducing disk patterns with small areas will be penalised, hence generating a repulsion effect. As for the Strauss process, the range parameter $r$ is considered known, and it cannot be taken directly into account because the likelihood $p(\yy|\theta)$ is not a differentiable function of it.\\

\noindent
In the following, we explore the structure of the galaxy pattern by assuming it is a realisation of an area-interaction process and by investigating its posterior distribution. Within this context, the posterior distribution gives important information related to the direction and the strength of points aggregation in the pattern. If the prior used is the uniform distribution, as in the previous situations, the maximum of the posterior is the maximum likelihood estimate. In this case, the posterior distribution gives, simultaneously, information about the general morphology of the pattern and about the model able to reproduce a pattern having sufficient statistics similar to the observed ones.\\

\noindent
For this purpose, the ABC Shadow algorithm was used in the following way. First, similarly to the summary statistics analysis, a finite set of values for the range parameters was fixed. From the galaxy point pattern, for each fixed value of $r$, the sufficient statistics of the area-interaction model were computed. Table~\ref{tableGalaxy} presents these values.\\

\begin{table}[!htbp]
\begin{center}
\begin{tabular}{|c|c|c|c|c|c|c|c|} 
\hline
\multicolumn{8}{|c|}{Data for the Galaxy pattern}\\
\hline
$r$            & 0.01 & 0.02 & 0.03 & 0.04 & 0.05 & 0.06 & 0.07\\
\hline
\multicolumn{8}{|c|}{$n(\yy) = 163$}\\
\hline
$- a_r(\yy)$   & 135.91 & 114.05 & 96.44 & 82.23 & 69.85 & 59.05 & 49.73\\
\hline
\end{tabular}
\end{center}
\caption{The observed sufficient statistics computed for the galaxy pattern, depending on the range parameter $r$. For all these parameters $n(\yy)$ remains constant, while $a_{r}(\yy)$ depends on $r$.}
\label{tableGalaxy}
\end{table}

\noindent
For each $r$ value an ABC Shadow algorithm was initialised with the computed sufficient statistics as the corresponding observed data. Except, the observed sufficient statistics, the parameters of all the algorithms were all the same. The prior density $p(\theta)$ was the uniform distribution on the interval $[2,12] \times [-5,5]$. Each time, the auxiliary variable was sampled using $250$ steps of a MH dynamics. The $\Delta$ and $n$ parameters were set to $(0.01,0.01)$ and $(100,100)$. The algorithm was run for $10^6$ iterations. Samples were kept every $10^3$ steps. This gave a total of $1000$ samples.\\

\begin{figure}[!htbp]
\begin{center}
\begin{tabular}{cc}
\epsfxsize=6cm \epsffile{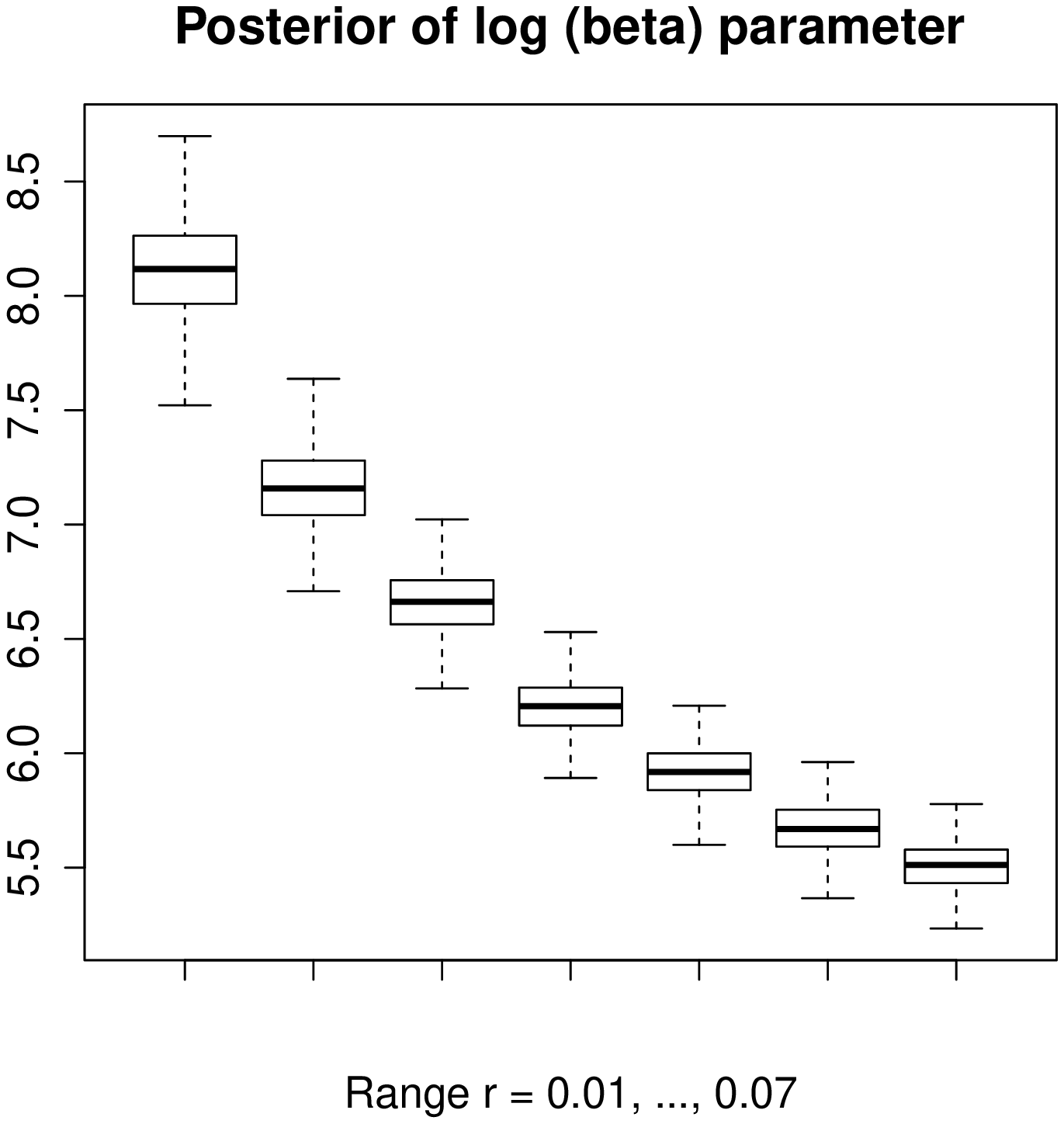} &
\epsfxsize=6cm \epsffile{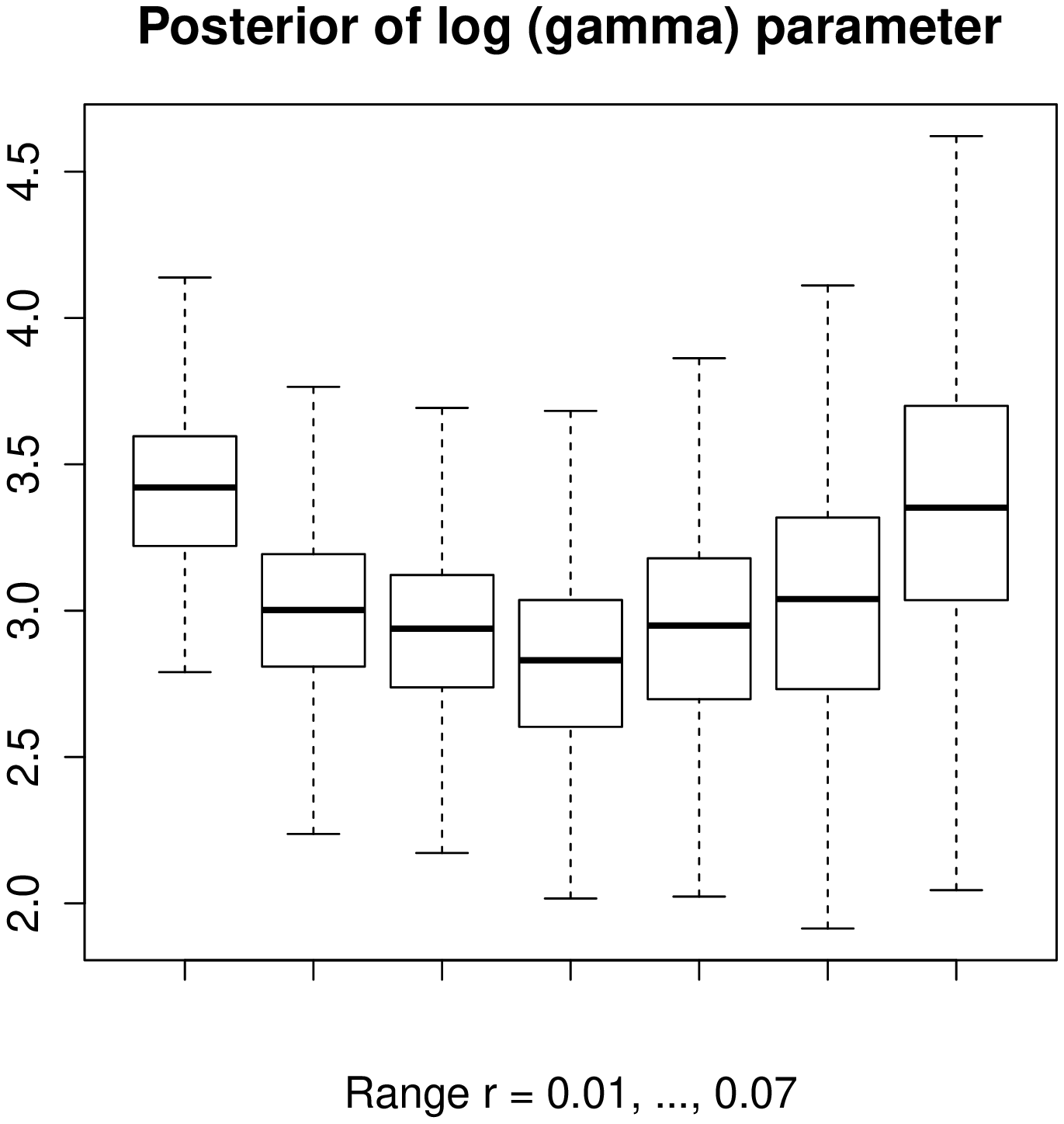}\\
\end{tabular}
\end{center}
\caption{Box-plots of the posterior distributions for the parameters of the area-interaction process estimated from the galaxy pattern, given different values for the interaction radius.}
\label{resultsGalaxy}
\end{figure}

\noindent
Figure~\ref{resultsGalaxy} shows the boxplots of the posterior distributions for the parameters $\log\beta$ and $\log\gamma$ conditionally on their corresponding range parameter. Since, the boxplots whiskers for the $\log\gamma$ posterior are far away from $0$ it may be concluded that there is a strong evidence of clustering of the pattern, for the considered ranges.

\noindent
The model parameters can be estimated for each radius. As in the previous examples, the MAP estimate equals the Maximum Likelihood (ML) estimate. In order to verify these results, the estimation errors are computed as in~\cite{LiesStoi03}. The results are shown in Table~\ref{tableAsymptotics}. For their computation, a simulation of $10000$ samples for each estimated model was performed. The asymptotic standard deviation approximates the difference between the maximum likelihood and the true model parameters, which are both unknown. The Monte Carlo standard deviation indicates here the difference between the MAP estimate and the maximum likelihood.\\
 
\begin{table}[!htbp]
\begin{center}
\begin{tabular}{|c|c|c|c|c|c|c|c|} 
\hline
\multicolumn{8}{|c|}{Estimation errors}\\
\hline
$r$            & 0.01 & 0.02 & 0.03 & 0.04 & 0.05 & 0.06 & 0.07\\
\hline
\multicolumn{8}{|c|}{Asymptotic standard deviation }\\
\hline
$\widehat{\sigma}_{\log\beta}$  & 0.20 & 0.17 & 0.13 & 0.11 & 0.10 & 0.10 & 0.08\\
$\widehat{\sigma}_{\log\gamma}$ & 0.26 & 0.28 & 0.27 & 0.30 & 0.34 & 0.40 & 0.52\\
\hline
\multicolumn{8}{|c|}{Monte Carlo standard deviation }\\
\hline
$\widehat{\sigma}_{\log\beta}^{MC}$  & 0.001 & 0.002 & 0.001 & 0.002 & 0.002 & 0.002 & 0.003\\
$\widehat{\sigma}_{\log\gamma}^{MC}$ & 0.002 & 0.003 & 0.004 & 0.006 & 0.008 & 0.012 & 0.024\\
\hline
\end{tabular}
\end{center}
\caption{Estimation errors computed for the MAP estimates obtained from the galaxy pattern.}
\label{tableAsymptotics}
\end{table}

\section{Conclusion and perspectives}
\noindent
The proposed method enlarges the panel of existing ABC techniques. To our best knowledge, it is the first time an ABC based statistical analysis is based on these point processes. The method was tested and tuned on known posteriors from the exponential family models. The results obtained on posteriors of point processes are coherent with the summary statistics analysis and the maximum likelihood estimation. Compared with the pseudo-likelihood based inference, our algorithm has the advantage that its output is a distribution that tends to approach the true posterior. The theoretical results hold for a class of models larger than the exponential family and they allow the use of different prior distributions. Clearly, the choice of the appropriate statistics and the construction of the test procedures for the algorithm within this general context, are open problems.\\

\noindent
As for all the other methods of sampling posteriors, the simulation of the auxiliary variable should be done exactly. Nevertheless, the numerical results obtained with less perfect samplers were satisfactory. The strong point of this method is that by getting close to the posterior distribution, no "useless" samples are produced. As perspectives, we mention range parameters estimation and model validation.\\

\section*{Aknowledgements}
This work was initiated during the stays of the first author at University Jaume I and INRA Avignon. The first author is grateful to D. Allard, Yu. Davydov, M. N. M. van Lieshout, J. M{\o}ller, E. Saar and the members of the Working Group ``Stochastic Geometry'' of the University of Lille, for useful comments and discussions. The work of the first author was partially supported by the GDR GEOSTO project. P. Gregori and J. Mateu were supported by grants P1-1B2012-52 and MTM2013-43917-P.\\

\section*{Appendix : proof of the results}
\subsection*{Proof of Theorem~\ref{first_theorem}}
\renewcommand{\labelenumi}{$($\roman{enumi}$)$}
\begin{enumerate}
\item 
Both density functions vanish outside the ball
$b(\theta, \Delta/2)$. For $\psi \in b(\theta, \Delta/2)$, the
integral mean value theorem applied to the denominator of $q_\Delta$
leads to
$$
q_\Delta(\theta \to \psi) = \frac{ \frac{f(\xx|\psi)}{c(\psi)} }{
V_\Delta \frac{f(\xx|\theta^*)}{c(\theta^*)} }
$$
for some $\theta^* \in b(\theta, \Delta/2)$. The positivity and the
continuity of the density $p$ (uniform continuity indeed, since
$\Theta$ is compact) allows us to do the following. Let $m(\xx) :=
\inf_{\phi \in \Theta} p(\xx | \phi) > 0$ since it is actually a minimum. For $A \in \cT_{\Theta}$ we have
\begin{align*}
\int_A & \left| q_{\Delta}( \theta \to \psi ) - U_{\Delta}( \theta
\to \psi ) \right| d \psi =
\int_{A \cap b(\theta, \Delta/2))} \left| \frac{ \frac{f(\xx|\psi)}{c(\psi)} }{ V_{\Delta} \frac{f(\xx|\theta^*)}{c(\theta^*)} } - \frac{1}{V_\Delta} \right| d \psi \leq \\
& \frac{1}{V_\Delta}\sup_{\phi \in \Theta} \frac{c(\phi)}{f(\xx|\phi)} \int_{A \cap b(\theta, \Delta/2)} \left| \frac{f(\xx|\psi)}{c(\psi)} - \frac{f(\xx|\theta^*)}{c(\theta^*)} \right| d \psi \leq \\
& \frac{\mu(A \cap b(\theta, \Delta/2))}{V_\Delta} m(\xx)^{-1} \sup_{d(\psi, \theta^*) < \Delta} \left| \frac{f(\xx|\psi)}{c(\psi)} - \frac{f(\xx|\theta^*)}{c(\theta^*)} \right| \leq \\
& m(\xx)^{-1} \sup_{d(\psi, \theta^*) < \Delta} \left|
\frac{f(\xx|\psi)}{c(\psi)} - \frac{f(\xx|\theta^*)}{c(\theta^*)}
\right|
\end{align*}
where the last supremum is independent of $\psi$ and $\theta^*$, and
approaches to $0$ as far as $\Delta$ does. With the regularity condition on $p(\xx|\cdot)$, we can tune up the inequality, using the (differential) mean value theorem, to
\begin{align*}
\int_A  \left| q_{\Delta}( \theta \to \psi ) - U_{\Delta}( \theta
\to \psi ) \right| d \psi & \leq m(\xx)^{-1} \Delta \sup_{\psi^* \in \Theta} \left\| D_{\Theta} p(\xx|\psi^*) \right\| \\ &:= C_1(\xx, p, \Theta) \Delta
\end{align*}
where $C_1(\xx, p, \Theta)$ is a constant depending on $\xx,p$ and $\Theta$.

\item 
As previously, the use of the integral mean value
theorem gives the result, since
\begin{align*}
\sup_{\psi \in \Theta} & \left| \frac{q_\Delta (\theta \rightarrow
\psi | \xx)}{q_\Delta (\psi \rightarrow \theta | \xx)} - \frac{
\frac{f(\xx|\psi)}{c(\psi)}\1_{b(\theta,\Delta/2)}(\psi) }{
\frac{f(\xx|\theta)}{c(\theta)}\1_{b(\psi,\Delta/2)}(\theta) }
\right| \leq \\
& \sup_{\psi \in b(\theta,\Delta/2)} \left| \frac{
\frac{f(\xx|\psi)}{c(\psi)} / (V_\Delta
\frac{f(\xx|\theta^*)}{c(\theta^*)})  }{
\frac{f(\xx|\theta)}{c(\theta)} / (V_\Delta
\frac{f(\xx|\psi^*)}{c(\psi^*)}) } - \frac{
\frac{f(\xx|\psi)}{c(\psi)} }{ \frac{f(\xx|\theta)}{c(\theta)} }
\right| \leq \\
& \sup_{\psi \in b(\theta,\Delta/2)} \left| \frac{
\frac{f(\xx|\psi)}{c(\psi)} }{ \frac{f(\xx|\theta)}{c(\theta)} }
\right| \left| \frac{ \frac{f(\xx|\psi^*)}{c(\psi^*)} }{
\frac{f(\xx|\theta^*)}{c(\theta^*)} } - 1 \right| \leq \\
& M(\xx) m(\xx)^{-2} \sup_{d(\theta^*,\psi^*) \leq \Delta} \left|
\frac{f(\xx|\psi^*)}{c(\psi^*)} -
\frac{f(\xx|\theta^*)}{c(\theta^*)} \right| \nonumber
\end{align*}
where $M(\xx) := \sup_{\phi \in \Theta} p(\xx | \theta) < \infty$ is
a maximum and $\theta^* \in b(\theta, \Delta/2)$, $\psi^* \in
b(\psi, \Delta/2)$ are obtained from the integral mean value
theorem. As in (i), under the regularity condition on $p(\xx|\cdot)$, this inequality can evolve to
\begin{eqnarray*}
\lefteqn{
\sup_{\psi \in \Theta} \left| \frac{q_\Delta (\theta \rightarrow
\psi | \xx)}{q_\Delta (\psi \rightarrow \theta | \xx)} - \frac{
\frac{f(\xx|\psi)}{c(\psi)}\1_{b(\theta,\Delta/2)}(\psi) }{
\frac{f(\xx|\theta)}{c(\theta)}\1_{b(\psi,\Delta/2)}(\theta) }
\right| \leq } \nonumber \\
& \leq & 
M(\xx) m(\xx)^{-2} \Delta \sup_{\psi^* \in \Theta} \left\| D_{\Theta} p(\xx|\psi^*) \right\| \\ 
& := & C_2(\xx, p, \Theta) \Delta \nonumber
\end{eqnarray*}
where $C_2(\xx, p, \Theta)$ is a constant depending on $\xx,p$ and $\Theta$.
\end{enumerate}

\subsection*{Proof of Proposition~\ref{first_proposition}}
\noindent
If $n = 1$, the definition of the transition kernels, the introduction of the term
$U_\Delta(\theta \to \psi) \alpha_i(\theta \to \psi) - U_\Delta(\theta \to
\psi) \alpha_i(\theta \to \psi)$, then the use of the triangle's inequality and
also the boundedness of functions $1_{A}(\cdot)$, $\alpha_i(\cdot)$ and
$\alpha_s(\cdot)$, allow us to write
\begin{align*}
| &P_i(\theta,A) - P_s(\theta,A) | \leq \\
   &\int_{\psi \in b(\theta, \Delta/2)} | q_\Delta(\theta \to \psi)
\alpha_i(\theta \to \psi) - U_\Delta(\theta \to \psi) \alpha_s(\theta 
\to \psi)
| d \psi \quad + \\
   & + 1_{A}(\theta) \int_{\psi \in b(\theta, \Delta/2)} | 
q_\Delta(\theta \to
\psi) [ 1 - \alpha_i(\theta \to \psi) ] - U_\Delta(\theta \to \psi) [ 1 -
\alpha_s(\theta \to \psi) ] | d \psi\\
   &\leq 3 \int_{\psi \in b(\theta, \Delta/2)} | q_\Delta(\theta \to \psi) -
U_\Delta(\theta \to \psi) |  d \psi +\\
   & + 2 \int_{\psi \in b(\theta, \Delta/2)} U_\Delta(\theta \to \psi) |
\alpha_s(\theta \to \psi) - \alpha_s(\theta \to \psi) | d \psi
\end{align*}

\noindent
Under the hypothesis of Theorem~\ref{first_theorem}, and then applying
Theorem~\ref{first_theorem}(i) and Corollary~\ref{first_corollary}, the
transition kernels of the ideal and the shadow Markov
chains, respectively are uniformly close as well (say for any
$\epsilon > 0$, there exists $\Delta(\epsilon,1) > 0$ such that we have $| P_s(\theta,A) - P_i(\theta,A) | < \epsilon$ if $\Delta \leq \Delta(\epsilon,1)$, independently of $\theta$ and $A$). \\

\noindent
If $p(\xx|\cdot) \in \cC^1(\Theta)$, then the previous inequalities can be completed to
\begin{equation*}
| P_i(\theta,A) - P_s(\theta,A) | \leq 3 C_1(\xx, p, \Theta) \Delta + 2 C_2(\xx, p, \Theta) \Delta  := C_3(\xx, p, \Theta) \Delta
\end{equation*}
giving a candidate expression for $\Delta_0(\epsilon, 1) := C_4(\xx, p, \Theta) \epsilon$. The constants $C_3$ and $C_4$ depend on $\xx,p$ and $\Theta$.\\

\noindent
For $n > 1$ we get by induction that
\begin{align*}
|P_i^{(n)}  - P_s^{(n)} | &\leq |P_i^{(n)}  - P_i^{(n-1)}  P_s^{(1)} | +
|P_i^{(n-1)}  P_s^{(1)}  - P_s^{(n)} | \\
&\leq |P_i^{(n-1)} | |P_i^{(1)}  - P_s^{(1)} | + |P_i^{(n-1)}  - 
P_s^{(n-1)} | |
P_s^{(1)} | \\
&\leq |P_i^{(1)}  - P_s^{(1)} | + |P_i^{(n-1)}  - P_s^{(n-1)} | \\
&\leq \cdots \leq n |P_i^{(1)}  - P_s^{(1)} | < \epsilon
\end{align*}
uniformly in $\theta$ and $A$ for all $\Delta$ such that
\begin{equation}
\Delta \leq \Delta_0(\epsilon, n) := \Delta_0(\epsilon/n, 1) = C_4(\xx, p, \Theta) \frac{\epsilon}{n}.
\label{delta_schedule}
\end{equation}

\bibliographystyle{plain}
\bibliography{abcpp}

\end{document}